\documentclass[12pt,reqno]{amsart}
\usepackage{amssymb,delarray}
\usepackage{amsfonts,amscd}
\usepackage{epsfig}
\usepackage[all]{xy}

\newcommand{\Sing}{\operatorname{Sing}}

\newtheorem{thm}{Theorem}
\newtheorem{lem}
{Lemma}
\newtheorem{prop}
{Proposition}
\newtheorem{claim}
{Claim}

{Corollary}
\newtheorem{rem}
{Remark}
{Question}
\newtheorem{ex-thm}{Theorem-Example}

{\catcode`\@=11
\gdef\n@te#1#2{\leavevmode\vadjust{%
 {\setbox\z@\hbox to\z@{\strut#1}%
  \setbox\z@\hbox{\raise\dp\strutbox\box\z@}\ht\z@=\z@\dp\z@=\z@%
  #2\box\z@}}}
\gdef\leftnote#1{\n@te{\hss#1\quad}{}}
\gdef\rightnote#1{\n@te{\quad\kern-\leftskip#1\hss}{\moveright\hsize}}
\gdef\?{\FN@\qumark}
\gdef\qumark{\ifx\next"\DN@"##1"{\leftnote{\rm##1}}\else
 \DN@{\leftnote{\rm??}}\fi{\rm??}\next@}}

\begin{document}
\baselineskip=13.7pt plus 2pt 

\title[]{On the variety of the inflection points of plane cubic curves}
\author[]{Vik.S. Kulikov}

\address{Steklov Mathematical Institute of Russian Academy of Sciences, Moscow, Russia}
 \email{kulikov@mi.ras.ru}

\dedicatory{} \subjclass{}

\keywords{}

\begin{abstract} In the paper, we investigate properties of the nine-dimensional variety of the inflection points of the plane cubic curves. The description of local monodromy groups of the set of inflection points near singular cubic curves is given. Also, it is given  a detailed description of the normalizations of the surfaces of the inflection points of  plane cubic curves belonging to general two-dimensional linear systems of cubic curves, The vanishing of the irregularity a smooth manifold birationally isomorphic to the variety of the inflection points of the plane cubic curves is proved.
\end{abstract}

\maketitle
\setcounter{tocdepth}{1}


\def\st{{\sf st}}


\setcounter{section}{-1}

\section{Introduction} \label{introduc}
The study of the properties of inflection points of plane nonsingular cubic curves has a rich and long history. Already in the middle of the 19th century, O. Hesse (\cite{He1}, \cite{He2}) proved that nine inflection points of a nonsingular plane cubic is a projectively rigid configuration of points, i.e. $9$-tuples of inflection points of the nonsingular plane cubic curves constitute an orbit with respect to the action of the group
$PGL (3,\mathbb C)$ on the set of $9$-tuples of the points of the projective plane. A subgroup $Hes$ in $PGL (3,\mathbb C)$ of order $216$ leaving invariant the set of the  inflection points of Fermat cubic, given in $\mathbb P^2$ by $z_1^3+z_2^3+z_3^3=0$, was defined by K. Jordan (\cite{J}) and named by him the {\it Hesse group}. The invariants of the group $Hes$  were described by G. Mashke in \cite{M}. A brief historical overview of the results relating to the inflection points of plane cubic curves can be found in \cite{A-D}.

Let  $F(\overline a,\overline z)=\displaystyle \sum_{0\leq i+j\leq 3}a_{i,j}z_1^iz_2^j\emph{}z_3^{3-i-j}$ be the homogeneous polynomial of degree three  in variables $z_1,z_2,z_3$ and of degree one in variables  $a_{i,j}$, $0\leq i+j\leq 3$. Denote by $\mathcal C\subset \mathbb P^{9}\times \mathbb P^2$ the complete family of plane cubic curves  given by equation $F(\overline a,\overline z)=0$.
Let $\kappa:\mathcal C\to  \mathbb P^{9}$ and $h:{\mathcal I}=\mathcal C\cap\mathcal H\to \mathbb P^{9}$ be the restrictions of the projection $\kappa_1: \mathbb P^{9}\times \mathbb P^2\to \mathbb P^{9}$ to $\mathcal C$ and $\mathcal I$, where
$$\displaystyle \mathcal H=\{ (\overline a,\overline z)\in \mathbb P^{9}\times \mathbb P^2\mid\det (\frac{\partial^2 F(\overline a,\overline z)}{\partial z_i\partial z_j})=0\} .$$

It is well known (see, for example, \cite{BK}), that for a generic point $\overline a_0\in \mathbb P^{9}$ the intersection of the curve $C_{\overline a_0}=\kappa^{-1}(\overline a_0)$ and its Hesse curve $H_{C_{\overline a_0}}$, given by \newline $\det (\frac{\partial^2 F(\overline a_0,\overline z)}{\partial z_i\partial z_j}) =0$, is the $9$ inflection points of $C_{\overline a_0}$. Therefore, $\deg h=9$.

Let $\mathcal B\subset \mathbb P^{9}$ be the discriminant hypersurface consisting of points $\overline a$ such that the curves $C_{\overline a}$ are singular. Then $h: {\mathcal I}\setminus h^{-1}(\mathcal B) \to \mathbb P^{9}\setminus \mathcal B$ is an unbranched covering of degree nine and, therefore, this covering determines a homomorphism
$h_ { * }: \pi_1(\mathbb P^{9}\setminus \mathcal B, \overline a_0) \to \mathbb S_{9}$ (here $\mathbb S_{9}$ is a symmetric group acting on the set $I_{\overline a}=C_{\overline a_0}\cap \mathcal I$). The group $\mathcal G=\text{Im}\, h_{*}$ is called {\it the monodromy group of the inflection points of plane cubic curves.} In \cite{Ha} has been proved

\begin{thm} \label{main} The group $\mathcal G$  is a group of order $216$. It is isomorphic to the group $Hes$ of projective transformations of the plane 
leaving invariant the set of  inflection points of  the Fermat curve of degree three.
\end{thm}

In this paper, we continue the study of the properties of the varieties of inflection points of plane curves 
begun in \cite{Ku-UMN}. In \S 1 we recall known facts about fundamental groups of complements to subvarieties of codimension one and on monodromy groups of dominant morphisms.  In \S 2 we give a description of the local monodromy groups of the set of inflection points of plane cubic curves near the points $\overline a\in \mathbb P^9$ parametrizing singular cubic curves (Propositions \ref{propx} -- \ref{loc3,2}). In \S 3 and \S 4, we give a detailed description of the normalizations of the surfaces of the inflection points of plane cubic curves belonging to generic two-dimensional linear systems of cubic curves, and also prove the main result of this paper (Theorem \ref{main2}), namely, we prove the vanishing of the irregularity of smooth varieties birationally isomorphic to the variety $\mathcal I$.

\section{The covering monodromy }
\subsection{On the fundamental groups of complements to subvarieties of codimen\-sion one.}
\label{f-loc}
Let $B\subset Z$ be a reduced closed codimension one subvariety of a simply connected smooth manifold $Z$, $\dim Z=K$. It is well known that the fundamental group $\pi _1(Z\setminus B, p)$ is generated by so-called {\it geometric generators} (or {\it bypasses} around $B$) $\gamma_q$, i.e., elements represented by loops $\Gamma_q$, $q\in B\setminus \Sing\, B$, of the following kind.
Let $L \subset Z$ be a germ of a smooth curve crossing transversely the variety $B$ at
$q\in B\setminus \Sing\, B$, and let $S_1 \subset L$ be a small radius circle centered at $q$. Right orientation of $Z$, defined by a complex structure,  determines the orientation on $S_1$ and in this case, $\Gamma_q$ is a loop consisting of the path $l$ lying in $Z\setminus B$ and connecting the point $p$ with some point $q_1\in S_1$, the loop $S_1$ (with the right orientation) starting and ending at the point $q_1$, and return  to the point $p$ along the path $l$. (Of course, the geometric generators $\gamma_q$ depend not only on the choice of the point $q$, but also on the choice of the path $l$; if we choose different paths and different points, the geometric generators $\gamma_q$ are conjugated to each other in $\pi _1(Z\setminus B,p)$ for points $q$ belonging to the same irreducible component of the variety $B$.)

In the case when $B$ is an hypersurface in $\mathbb P^k$ and $\mathbb P^n$ is a $n$-plane in $\mathbb P^k$, which is in general position with respect to $B$, then Zariski -- van Kampen Theorem states  that the homomorphism $i_*:\pi_1(\mathbb P^n\setminus B)\to \pi_1(\mathbb P^k\setminus B)$, induced by inclusion $i:\mathbb P^n\hookrightarrow \mathbb P^k$ is an isomorphism when $n\geq 2$ and epimorphism when $n=1$. In the case $n=1$, the group $\pi_1(\mathbb P^1\setminus B,p)$ is generated by $d=\deg B$ bypasses $\gamma_1,\dots, \gamma_d$ in $\mathbb P^1\setminus B$ around the points belonging to the set $\{ q_1,\dots, q_d\}=B\cap \mathbb P^1$. An ordered set $\{ i_*(\gamma_1),\dots, i_*(\gamma_d)\}$ is called {\it a good geometric base} of the group
$\pi_1(\mathbb P^k\setminus B,p)$ if the product $\gamma_1\cdot .\, .\, .\, \cdot \gamma_d$ is equal to the neutral element in
$\pi_1 (\mathbb P^1\setminus b,p)$.

Let $m=m_o(B)$ be the {\it multiplicity of singularity} of the hypersurface $B$ at the point $o$ and let $z_1,\dots ,z_k$ be local analytic coordinates in a certain analytic neighborhood $V\subset Z$ of the point $o=(0,\dots,0)\in B$ chosen so that the  intersection number $(L_0,B)_o$ at the point $o$ of the "line"\, $L_0\subset V$, given by $z_1=\dots=z_{k-1}=0$, and $B$ is equal to $m$.

Let $$\Delta^k_{\overline{\delta}} =\{ (z_1,\dots,z_k)\in V\mid |z_i|< \delta_i\,\, \text{for}\, \, i=1,\dots, k \} \subset Z$$ be a polidisk of multiradius $\overline{\delta}=(\delta_1,\dots,\delta_k)>0$ centered at the point $o$. Each time, we assume that multiradius $\overline{\delta}=(\delta_1,\dots,\delta_k)$ is chosen so that the restriction of the projection $\text{pr}:(z_1,\dots,z_{k-1},z_k)\mapsto (z_1,\dots,z_{k-1})$ to $\Delta_{\overline{\delta}}^k\cap B$ is a proper finite map of degree $m$. Let a "line"\, $L\subset {\Delta}_{\overline{\delta}}^k$, defined by equations $z_1=c_1<\delta_1, \dots, z_{k-1}=c_{k-1}<\delta_{k-1}$ intersect $B$ in $m$ different points $q_1,\dots, q_m$ and let $p=(c_1,\dots, c_{k-1},c_k)\in L\setminus B$. The fundamental group $\pi_1 (L\setminus B,p)$ is the free group of rank $m$ and it is generated by $m$ bypasses $\gamma_1,\dots,\gamma_m$ around the points $q_1\dots, q_m$. An ordered set
$ \{\gamma_1,\dots, \gamma_m\}$ is called a {\it  good geometric base} of $\pi_1(L\setminus B, p)$ if the product
$\gamma_1\cdot .\, .\, .\, \cdot \gamma_m$ is equal to the element in $\pi_1 (L\setminus B,p)$ represented by the circuit along the circle $\Delta= \{|z_k| =\delta_k\}$. The ordered set $\{ i_*(\gamma_1),\dots, i_*(\gamma_m)\}$ is also called  a {\it  good geometric base} of the group
$\pi_1 (\Delta^k_{\overline{\delta}}\setminus B,p)$. It is easy to show that the elements of a good geometric base generate the group
$\pi_1 (\Delta^k_{\overline{\delta}}\setminus B, p)$).

For multiradii $\overline{\delta}_1$ and $\overline{\delta}_2$, $\overline{\delta}_1\leq \overline{\delta}_2$, the  imbedding  $\Delta_{\overline{\delta}_1}^k\subset \Delta_{\overline{\delta}_2^k}$ induces a homo\-morphism $i_*:\pi_1(\Delta^k_{\overline{\delta}_1}\setminus B)\to\pi_1(\Delta^k_{\overline{\delta}_2}\setminus B)$ of the fundamental groups.
The following theorem is well known (see, for example, \cite{HLe}).

\begin{thm} \label{loc}
There is such a multiradius $\overline r>0$ that for any multiradius $\overline{\delta}\leq \overline r$ the homomorphism $i_*:\pi_1(\Delta^k_{\overline{\delta}}\setminus B)\to \pi_1(\Delta^k_{\overline r}\setminus B)$, induced by the polydisk embedding is an isomorphism.
\end{thm}
The group $\pi_1^{loc} (B,o):=\pi_1 (\Delta^k_{\overline r}\setminus b)$ is called {\it local fundamental group} of $B$ at $o$.

Let $\Pi$ be a linear subspace of $\Delta_{\overline{\delta}}^k$, $\dim \Pi=n$, passing through $o$ and being at that point in general position with respect to $b$. The following theorem is a direct consequence of the results of \cite{H-le1}.
\begin{thm} \label{loc1}
The homomorphism $i_*:\pi_1^{loc} (B\cap \Pi,o)\to \pi_1^{loc} (B,o)$, induced by an embedding $i:\Pi\hookrightarrow \Delta^k_{\overline{\delta}}$, is an isomorphism for $n\geq 3$ and an epimorphism for $n=2$.
\end{thm}

The following claims are well known.
\begin{claim} \label{Exam-c}  Let the curve $B\subset \Delta^2_r$ is given by  $z_1^3+z_2^2=0$, i.e., $B$ has a singularity of type $A_{2}$ at the point $o$. Then the group $\pi_1^{loc} (B,o)$ has the following presentation:
\begin{equation} \label{a2n-1}
\pi_1^{loc} (B,o)=\langle \gamma_1,\gamma_2\ : \gamma_1\gamma_2\gamma_1=\gamma_1\gamma_2\gamma_1\rangle, \end{equation}
where $ \{ \gamma_1,\gamma_2\}$ is a good geometric base of $\pi_1(\Delta_r^2\setminus B)$. \end{claim}

\begin{claim} \label{ex2} Let $B\subset \Delta^k_r$ be given by equation $z_1\dots z_n=0$, $n\leq k$, i.e., $B$ is a divisor with normal crossings at $o$. Then the group $\pi_1^{loc} (B,o)\simeq \mathbb Z^m$ is a free abelian group and it is generated by bypasses $\gamma_i$, $i=1,\dots, m$, around hypersurfaces  $\{ z_i=0\}$. \end{claim}

\begin{lem} \label{sigma} {\rm (}\cite{Ku-S}{\rm )} Let $(B,o)$ be a germ of $B$ in $\Delta_r^2$ and let
$\sigma:\widetilde{\Delta_r^2}\to \Delta^2_r$ be $\sigma$-process with center at $o$, and $E=\sigma^{-1}(o)$ its exceptional divisor. Then the  conjugacy class in the group $\pi_1(\widetilde{\Delta_r^2}\setminus \sigma^{-1} (B),\sigma^{-1} (p))=\pi_1 (\Delta_r^2\setminus B,p)$ of the bypass $\gamma$ around $E$ contains the element $\Delta=\gamma_1\cdot \, . \, . \, .\, \cdot\gamma_m$, where
$\{ \gamma_1,\dots,\gamma_m\}$ is a good geometric base of the group $\pi_1(\Delta_r^2\setminus B,p)$.
\end{lem}

\subsection{On the monodromy of dominant morphisms.} \label{monodr}
Let $f:X\to Z$ be a dominant proper holomorphic map  of a complex-analytic variety $X$ on smooth manifold $Z$, $\dim X=\dim Z=k$. In this case, there is a subvariety $B\subset Z$ of codimension one (which we call the {\it discriminant} of  $f$) such that $f: Y=X\setminus f^{-1}(B)\to Z\setminus B$ is a finite unramified covering. Note that $Y$ is a smooth manifold.

Let the degree of  $f:Y\to Z\setminus B$ be $n$. The map $f$ defines a homomorphism $f_*:\pi_1(Z\setminus B,p)\to \mathbb S_n$ (called the {\it monodromy} of $f$) whose image $G_f:=f_*(\pi_1(Z\setminus B,p))$  is called the {\it monodromy group} of  $f$ and it is a subgroup of the symmetric group $\mathbb S_n$ acting on the fibre $f^{-1}(p)=\{ p_1,\dots,p_n\}$ in the following way. A loop $\Gamma\subset Z\setminus B$ representing an element $\gamma \in\pi_1(Z\setminus B,p)$ can be lifted to $Y$ and as a result,  we get $n$ paths $\Gamma_1,\dots,\Gamma_n\subset Y$ starting and ending at the points of the fibre $f^{-1}(p)$. Therefore, this lift defines the action of $f_*(\gamma)$ on $f^{-1}(p)$, which for each $i=1,\dots, n$ maps the starting point $p_i$ of the path $\Gamma_i$ to the end point of that path.

By Grauert - Remmert - Riemann - Stein Theorem (\cite{St}) each epimorphism $\varphi_*:\pi_1(Z\setminus B)\to G\subset \mathbb S_n$ uniquely determines a finite covering $\varphi:X\to Z$ (which we call the {\it Stein covering} associated to $ \varphi_*$) of degree $\deg \varphi =n$, unramified over $Z\setminus B$, where $X$ is a normal variety, and the group $G\subset \mathbb S_n$ is the monodromy of  $\varphi$, and if $\varphi_*=f_*$ and $Y$ is a normal variety, then there is a dominant map
$\psi :Y\to X$ (the factor map contracting the connected components of the fibres of $f$ at points of $X$) such that $f=\varphi\circ\psi$. In the case when $n=|G|$ and the embedding $G\subset \mathbb S_ {|G|}$ is the Cayley embedding, i.e., $G$ acts on itself by multiplication on the right, the group $G$ acts on $X$ such that $X/G=Z$ and in this case  $f$ is the {\it Galois covering}. The following proposition is true.

\begin{prop} \label{gal-st} If $G\subset \mathbb S_n$ is the monodromy group of Stein covering $\varphi :X\to Z$
and $\mathbb S_n$ acts on the set $\{q_1,\dots,q_n\}$, then for the Galois covering $\widetilde{\varphi}: Y\to Z$ with the Galois group $G$ there is a Galois covering $\overline{\varphi}:Y\to X$ with Galois group $G^1=\{g\in G\mid g(q_1)=q_1\}$ such that
$\widetilde{\varphi}=\varphi\circ \overline{\varphi}$.
\end{prop}

In the notation used in subsection \ref{f-loc}, let $o$ be a point in a hypersurface $B$. The choice of the path connecting the base points of the fundamental groups determines  homomorphisms $i_*:\pi_1^{loc} (B,o) \rightarrow \pi_1(Z\setminus B)$ and  $f_{*,loc}=f_*\circ i_*:\pi_1^{loc} (B,o)\to G$ defined uniquely up to conjugation in the group $G$.
The group $G_{o}:=Im\, f_{*, loc}\subset G\subset \mathbb S_n$ is called the {\it local monodromy group} of $f$ at $o$.

The following claim directly follows from the proof of Grauert - Remmert - Riemann - Stein Theorem.
\begin{claim} \label{connect} Let $J_1,\dots, J_m$ be orbits of the action of the local monodromy group $G_o\subset \mathbb S_n$ of Stein covering $\varphi$ on the set $\{1,\dots,n\}$, where $n_j$ is the cardinality  of the orbit $J_j$, $n_1 + \dots +n_m=n$.  Then the preimage $\varphi^{-1} (V)$ of a sufficiently small neighborhood $V\subset Z$ of the point $o$ splits into a disjoint union of $m$  irreducible normal neighborhoods $U_j$. The degree of restriction of covering $\varphi$ to $U_j$ at $q_j=U_j\cap \varphi^{-1} (o)$ is equal to $n_j$. In particular, if $n_j=1$, then $U_j$ is nonsingular at $q_j$.

If $\varphi$ is a Galois cover and if the neighborhood of $V$ is small enough, then $\varphi^{-1}(V)=\bigsqcup_{j=1}^{i_{o}}U_j$ is a disjoint union of $i_{o}$ connected normal pairwise biholomorphic surfaces {\rm (}possibly having a singular  point $q_j\in U_j$), where $i_{o}=(G:G_o)$ is the index of the subgroups $G_{o}$ in $G$. The restriction of covering $\varphi$ for each neighborhood $U_j$ is a Galois covering with Galois group conjugated with the group $G_o$.
\end{claim}

Below we will need the following
\begin{claim}\label{deg3} Let $f:U\to \Delta_r^2$ be a proper dominant holomorphic map of degree three of a smooth surface $U$ to the bidisk $\Delta_r^2$ ramified in coordinate axes $B=\{ z_1z_2=0\}$ and whose fibre $L=f^{-1}(o)\simeq \mathbb P^1$ is a nonsingular rational curve. Assume that the images $f_*(\gamma_1)$ and $f_*(\gamma_2)$ of geometric generators $\gamma_1$ and $\gamma_2$ of the group $\pi_1(\Delta^2_r\setminus B)$ in the monodromy group $G\subset \mathbb S_3$ of  $f$  are mutually inverse elements in $G\simeq \mathbb Z_3$. Then the self-intersection number  of the curve $L$ in $U$ is  $(L^2)_U=-3$.
\end{claim}
\proof Let $\varphi: V\to \Delta_r^2$ be the Stein covering associated with epimorphism $f_*$, $\deg \varphi =3$. Then the factorization map $\psi: U\to V$ is the resolution of the singular point $o'=\varphi^{-1}(o)$ of $V$, where $o=(0,0)$ is the origin in $\Delta_r^2$. On the other hand, let $\sigma: \widetilde{\Delta}^2_r\to \Delta_r^2$ be the $\sigma$-process centered at $o$, $E=\sigma^{-1}(o)$. Let $\phi:W\to\widetilde{\Delta}^2_r$ be the Stein covering associated with epimorphism
$f_*:\pi_1(\Delta^2_r\setminus B)=\pi_1(\widetilde{\Delta}^2_r\setminus (\sigma^{-1}(B)\cup E))\to \mathbb S_3$. The covering $\phi$ is ramified with multiplicity three only along the smooth curve $\sigma^{-1} (B)$, the proper preimage of the curve $B$, since $f_*(\gamma)=Id$ according to Lemma \ref{sigma}, where $\gamma$ is a bypass around the curve $E$. Hence, $W$ is a nonsingular manifold and it is easy to see that $L_1=\phi^{-1}(E)$ is a smooth irreducible curve, $(L_1^2)_W=\deg \phi \cdot (E^2)_{\widetilde{\Delta}_r^2}=-3$.
According to Grauert -- Remmert -- -- Riemann -- Stein Theorem,  for $\sigma\cdot\phi: W\to \Delta_r^2$ there is a canonical factorization  $\psi_1: W\to V$,
also being the resolution of the singular point $o'\in V$. Therefore Claim \ref{deg3} follows from the uniqueness of the minimal resolution of singular points of two-dimensional complex varieties.   \qed

\section{Local monodromy groups of morthism $h$} \label{mon-loc}
\subsection{On the inflection points of plane cubic curves.}
Consider a cubic curve $C_{\overline a_0}\subset \mathbb P^2$ given by equation $F(\overline a_0,\overline z)=0$ and its Hessian curve $H_{C_{\overline a_0}}\subset \mathbb P^2$ given by equation
$\det (\frac{\partial^2 F(\overline a_0,\overline z)}{\partial z_i\partial z_j})=0$. By definition, the {\it inflection points} of  $C_{\overline a_0} $ are the points belonging to the intersection  $C_{\overline a_0}\cap H_{C_{\overline a_0}}$.
\begin{claim}\label{sing-inf} The singular points of a cubic curve $C_{\overline a_o}$ are its inflection points.
\end{claim}
\proof If $s=(0,0,1)$ is a singular point of a cubic curve $C_{\overline a_0}$, then $C_{\overline a_0}$ is given by equation of the form $z_3\sum_{i+j=2} c_{i,j}z_1^iz_2^j+\sum_{i+j=3}c_{i,j}z_1^iz_2^j=0$. It is easy to see that the elements of the third column of the Hesse matrix on the left side of this equation vanish at $s$.  \qed
\begin{claim}\label{lint-inf} The points of a line $L$ being a component of a cubic curve $C_{\overline a_o}$ are inflection points of  $C_{\overline a_o}$.
\end{claim}
\proof Let the line $L$ be given by equation $z_1=0$. In this case, the cubic curve $C_{\overline a_0}$ is given by  equation of the form $z_1\sum_{0\leq i+j\leq 2}c_{i, j}z_1^iz_2^jz_3^{2-i-j}=0$. It is easy to see that only the elements of the first column and the first row of the Hesse matrix of the left side of this equation can have non-zero values at the points of $L$. Therefore, Hessian vanishes at the points of $L$.  \qed

\subsection{Equisingular stratification of the variety of singular cubic curves.}\label{strat}
Denote $\mathcal S:=\{ (\overline a,\overline z)\in \mathbb P^9\times\mathbb P^2\mid \overline z\in Sing\, C_{\overline a}\}$. We have $\text{pr}_1(\mathcal S)=\mathcal B$. The discriminant $\mathcal B\subset \mathbb P^9$ is irreducible \cite{H-P} and it has the natural stratification:

\begin{picture}(300,80)
\put(57,40){${\mathcal B}_7$}
\put(72,45){\vector(1,0){30}}
\put(107,40){${\mathcal B}_5$}
\put(123,45){\vector(1,0){30}}
\put(172,47){\vector(2,1){30}} \put(172,43){\vector(2,-1){30}}
\put(207,60){$\mathcal B_{3,1}$}
\put(207,20){$\mathcal B_{3,2}$}\put(222,30){\vector(1,1){30}}
\put(158,40){$\mathcal B_4$}
\put(226,65){\vector(1,0){28}}
\put(226,25){\vector(1,0){28}}
\put(257,60){$\mathcal B_{2,1}$}
\put(257,20){$\mathcal B_{2,2}$}
\put(278,27){\vector(2,1){30}}
\put(278,63){\vector(2,-1){30}}
\put(315,40){$\mathcal B_1\subset \mathcal B,$}
\end{picture}
\newline where $\mathcal B_7$ is the variety parametrizing the  cubic curves consisting of  triple lines; $\mathcal B_5$ is the variety parametrizing the cubic curves consisting of two lines, one of which is included with multiplicity two; $\mathcal B_4$ is the variety parametrizing the cubic curves consisting of three lines with a common point; $\mathcal B_{3,1}$ is the variety parametrizing the cubic curves consisting of three lines  in general position; $\mathcal B_{3,2}$ is the variety parametrizing the cubic curves consisting of smooth conics and lines touching each other; $\mathcal B_{2,1}$ is the variety parametrizing the cubic curves consisting of transversely intersecting conics and lines; $\mathcal B_{2,2}$ is the variety parametrizing the cuspidal rational cubic curves; $\mathcal B_1$ is the variety parametrizing the nodal rational cubic curves; arrows denote the adjoint of the strata and for the strata $\mathcal B_i$ and $\mathcal B_{i,j}$ the index $i$ is equal to the codimension of these strata in $\mathbb P^9$.
Note that all strata are smooth manifolds and, moreover, these strata are orbits of the action of the group $\text{PGL}(3,\mathbb C)$ on the space of plane cubic curves. Therefore, $\mathcal B$ is an irreducible variety and therefore all geometric generators of the groups
$\pi_1(\mathbb P^9\setminus \mathcal B)$ belong to the same conjugacy class of elements of this group and generate it. In addition, we have the following
\begin{claim}\label{equi-str}
\begin{itemize}
\item[$(i)$] For any two points $\overline a_1$ and $\overline a_2$ belonging to the same stratum $\mathcal B'$ of the equisingular stratification of variety $\mathcal B$, the local monodromy groups $\mathcal G_{\overline a_1}$ and
$\mathcal G_{\overline a_2}$ are conjugated in $\mathcal G$.
\item[$(ii)$] If a stratum  $\mathcal B'$ adjoints to a stratum $\mathcal B''$, then the local monodromy group $\mathcal G_{\overline a_1}$ at $\overline a_1$ contains a subgroup conjugated in $\mathcal G$ with the local monodromy groups $\mathcal G_{\overline a_2}$ at the points $\overline a_2\in \mathcal B^{\prime\prime}$.
\end{itemize}
\end{claim}
\proof The group $\text{PGL}(3,\mathbb C)$ is a connected complex manifold. Therefore claim $(i)$ follows from that each stratum $\mathcal B'$ is an orbit under the action of $\text{PGL}(3,\mathbb C)$ on the variety of plane cubic curves.

Claim $(ii)$ follows from that for any sufficiently small neighbourhood $V_1\subset \mathbb P^9$ of a point $\overline a_1\in \mathcal B'$ there is a point $\overline a\in\mathcal B^{\prime\prime}$ lying in $V_1$ together with some its neighbourhood $V$. \qed \\

The space $\mathbb P^9$ has a natural covering by ten affine spaces $\mathbb C^9_{i,j}=\{ \overline a\in \mathbb P^9 \mid a_{i,j}\neq 0\}$, $0\leq i+j\leq 3$. The coordinates of the points
$\overline a$ in each space $\mathbb C^9_{i_0,j_0}$ are the coefficients in equations $$\sum_{0\leq i+j\leq 3}a_{i,j}z_1^iz_2^jz_3^{3-i-j}=0$$
of plane cubic curves normalized by condition $a_{i_0,j_0}=1$.

Properties of the discriminant $\mathcal B$ at the points corresponding the nodal cubic curves were investigated in \cite{Se}.

\begin{prop} \label{tangent} {\rm (\cite{Se})} Let $C_{\overline a_0}$ be a nodal cubic curve. Then  $\mathcal B$ is a divisor with normal crossings in some analytic neighbourhood of the point $\overline a_0$.
\end{prop}

Below,  to calculate the local monodromy groups of morphism $h$ at the points corresponding to the nodal cubic curves, we need  more detail description of the discriminant $\mathcal B$ at these points. For this, consider a nodal cubic curve $C_{\overline a_0}$ given in non-homogeneous coordinates $x=\frac{z_1}{z_3}$, $y=\frac{z_2}{z_3}$ in $\mathbb C^2\subset \mathbb P^2$ by equation
$$ xy+\sum_{i+j=3}c_{i,j}x^iy^j=0.$$
The point $\overline a_0$ lies in the affine space $\mathbb C^9_{1,1}$ and its coordinates $a_{0,0}$, $a_{1,0}$, and $a_{0,1}$
vanish. 
Obviously, the variety $\mathcal S\cap (\mathbb C^9_{1,1}\times \mathbb C^2)$ is given by equations
\begin{equation} \label{node} \begin{array}{l} a_{0,0}+a_{1,0}x+a_{0,1}y +xy +a_{2,0}x^2+a_{0,2}y^2+\sum_{i+j=3}a_{i,j}x^iy^j=0 \\
a_{1,0}+y+2a_{2,0}x+\sum_{i+j=3}ia_{i,j}x^{i-1}y^j=0 \\
a_{0,1}+x+2a_{0,2}y+\sum_{i+j=3}ja_{i,j}x^{i}y^{j-1}=0
\end{array}\end{equation}
It follows from the type of equations (\ref{node}) that the eight-dimensional variety $\mathcal S$ is smooth at the point $(\overline a_0,\overline z_0)$, where $\overline z_0=(0,0,1)$ and for $2\leq i+j\leq 3$, $(i,j)\neq (1,1)$ the functions $a_{1,0}$, $a_{0,1}$, and $a_{i,j}-c_{i,j}$ are local parameters in $\mathcal S$ at $(\overline a_0,\overline z_0)$. Therefore the image $\text{pr}_1(S\cap V)\subset \mathbb C^9_{1,1}$ is non-singular at $\overline a_0$, where $V$ is some sufficiently small analytic neighbourhood in $\mathbb C^9_{1,1}\times \mathbb C^2$ of the point $(\overline a_0,\overline z_0)$. In addition, it is easy to see that $a_{0,0}=0$ is an equation of the tangent space of $\text{pr}_1(S\cap V)$ at  $\overline a_0$.

The strata  $\mathcal B_1$, $\mathcal B_{2,1}$, and $\mathcal B_{3,1}$ of the discriminant hypersurface $\mathcal B$ parametrize the nodal cubic curves. These strata are orbits in $\mathbb P^9$ under the action of the group $PGL(3,\mathbb C)$. Therefore, without loss of generality, we can assume that  a nodal cubic $C_{\overline a_0}$ is given by one of the following equations:
$$ \begin{array}{ll} 1) & z_1z_2z_3+z_1^3+z_2^3=0  \,\, \text{if}\,\, \overline a_0\in \mathcal B_1, \\
2) & z_1z_2z_3+z_1^3=0\quad \quad \,\, \, \text{if}\,\, \overline a_0\in \mathcal B_{2,1}, \\
3) & z_1z_2z_3=0\quad \quad \quad \quad \,\,\,\, \text{if}\,\, \overline a_0\in \mathcal B_{3,1}.
\end{array} $$
In all three cases, the point $\overline a_0$ belongs to the affine space $\mathbb C^9_{1,1}$ and it follows from the above that
in case 1) $\mathcal B$ is a nonsingular hypersurface at $\overline a_0$; in case 2) $\mathcal B$ locally is a union of two nonsingular at point $\overline a_0$ hypersurfaces intersecting transversally at that point (tangent spaces to these hypersurfaces are given by equations $a_{0,0}=0$ and $a_{0,3}=0$); in  case 3) $\mathcal B$ is locally the union of three non-singular  hypersurfaces at  $\overline a_0$ intersecting transversally at this point (the tangent space to these hypersurfaces is defined by the equations $a_{0,0}=0$, $a_{0,3}=0$, and $a_{3,0}=0$).

\subsection{Properties of the monodromy group $\mathcal G=Hes$ of morphism $h$.}\label{Hes}
The group ${Hes}$ was originally defined as a subgroup of the group $PGL(3,\mathbb C)$ of linear transformations of $\mathbb P^2$ leaving  invariant the {\it Hesse pencil}, i.e., a one-dimensional linear system of plane cubic curves defined by  equation
\begin{equation}\label{Fer3} C_{(t_1,t_2)}: \quad t_1(z_1^3+z_2^3+z_3^3)+ t_2z_1z_2z_3=0, \quad (t_1,t_2)\in \mathbb P^1 ,\end{equation}
The Hesse pencil has nine fixed points
$$\begin{array}{lll}
q_1=(0,1,-1),\,\, & q_4=(0,1,-\omega),\,\, & q_7=(0,1,-\omega^2), \\
q_2=(1,0,-1),\,\, & q_5=(1,0,-\omega^2),\,\, & q_8=(1,0,-\omega), \\
q_3=(1,-1,0),\,\, & q_6=(1,-\omega,0),\,\, & q_9=(1,-\omega^2,0), \end{array} $$
where $\omega=e^{2\pi i/3}$ is a primitive cubic root of unity, and a direct checking shows that these nine points are inflection points of each of the nonsingular curves entering the Hesse pencil. Therefore, the group ${Hes}\subset PGL(3,\mathbb C)$ can be defined as a subgroup of projective transformations leaving invariant the set of inflection points of the Fermat cubic $F=C_{(1,0)}$.

The properties of the group $Hes$ are well known (see, for example, \cite{CM}). Its order is equal to $216$ and the action of the group ${Hes}$ on the nine inflection points of the curve $F$ defines an embedding ${Hes}\subset \mathbb S_9$ such that ${Hes}$ is a $2$-transitive subgroup of the symmetric group $\mathbb S_9$ generated by permutations $g_0=(1,2,4)(5,6,8)(3,9,7)$ and $g_1=(4,5,6)(7,9,8)$. The subgroups
$Hes^i=\{ g\in Hes \mid g(q_i)=q_i\}$, $i=1,\dots,9$, consisting of elements of $Hes$  leaving fixed the point $q_i$, are conjugated to each other in the group $Hes$ and they are isomorphic to the group $SL(2,\mathbb Z_3)$. Their order is $24$. As a subgroup in
$\mathbb S_9$, the group $Hes^1$ is generated by permutations $g_1$ and $g_2=g_0g_1g_0^{-1}=(2,8,5)(3,6,9)$ and as an abstract group, it has the following presentation
$$ Hes^1= \{ g_1,g_2 \mid g_1g_2g_1=g_2g_1g_2,\,\, g_1^3=1\}. $$
For further we note that the class in $Hes^1$ of elements conjugated to $g_1$ consists of four elements:
$g_1, g_2,g_1g_2, g_1^{-1}, g_2g_1g_2^{-1}$, and any pair of these four elements generates the group $Hes^1$ and satisfies the ralation $xyx=yxy$. Also note (see \cite{S-T}) that the group $Hes^1$ has an exact linear presentation in $GL(2,\mathbb C)$ (group no. 4 in \cite{S-T}, see also \cite{Ku3}), in which images of $g_1$ and $g_2$ are reflections.

The group $Hes^1$ is naturally embedded in the symmetric group
$\mathbb S_8$ acting on $\{ q_2,\dots, q_8\}$.  We number (left) co-sets in $Hes^{1}$ of its subgroup $\langle g_1\rangle$ generated by $g_1$ as follows:
$$\begin{array}{llll}  c_2=\langle g_1\rangle, & c_3= g^2_2g_1g_2^2\langle g_1\rangle, & c_4=g_1^2g^2_2\langle g_1\rangle, & c_5=g_2^2\langle g_1\rangle,\\ c_6=g_1g_2^2\langle g_1\rangle, & c_7= g_1g_2\langle g_1\rangle, & c_8= g_2\langle g_1\rangle, & c_9= g_1^2g_2\langle g_1\rangle.\end{array}$$
The group $Hes^1$ acts on the set $\{ c_2,\dots,c_9\}$.

Denote $Hes^{1,2}=\{ g\in Hes^1\subset \mathbb S_8\mid g(q_2)=q_2\}$.
\begin{claim} \label{norm-cusp}  We have $Hes^{1,2}=\langle g_1 \rangle$ and the action of the group
$$Hes^1=\langle g_1=(4,5,6)(7,9,8), g_2=(2,8,5)(3,6,9)\rangle\subset \mathbb S_8$$ on the set $\{ c_2,\dots,c_9\}$ coincides with the action on the set $\{ q_2,\dots,q_9\}$.
\end{claim}
\proof Writing all $24$ permutations belonging to the group $Hes^1\subset \mathbb S_8$, it is easy to verify the validity of this statement. \qed

\subsection{On the local monodromy groups of the morphism $h$ at points corres\-pond\-ing to nodal cubic curves.}

\begin{prop} \label{propx} $(i)$ The local monodromy groups 
of morphism $h$ at the points
$\overline a_{1}\in \mathcal B_{1}$ are abelian groups isomorphic to the group $\mathbb Z_3$ and they are conjugated in $Hes$ to the subgroup
$\mathcal G_{{1}}\subset \mathbb S_{9}$ generated by permutation
$g_2=(2,8,5)(3,6,9)$.

$(ii)$ The local monodromy groups 
of morphism $h$ at the points $\overline a_{2,1}\in \mathcal B_{2,1}$ are abelian groups isomorphic to the group $\mathbb Z_3^2$ and they are conjugated in $Hes$ to the subgroup $\mathcal G_{{2,1}}\subset \mathbb S_{9}$ generated by permutations
$g_2=(2,8,5)(3,6,9)$ and $g_3=(1,4,7)(3,9,6)$.

$(iii)$ The local monodromy groups of morphism $h$ at points $\overline a_{3,1}\in \mathcal B_{3,1}$ are abelian  groups isomorphic to the group $\mathbb Z_3^2$ and they are conjugated in $Hes$ to the subgroup $\mathcal G_{{3,1}}\subset \mathbb S_{9}$ generated by permutations
$g_2=(2,8,5)(3,6,9)$ and $g_3=(1,4,7)(3,9,6)$, i.e., coincide with the group $\mathcal G_{2,1}$.
\end{prop}
\proof We start the proof with statement $(iii)$. Without loss of generality, we can assume that a nodal curve $C_{\overline a_{3,1}}$ is given by equation $z_1z_2z_3=0$.

Consider a three-dimensional family $\mathcal C_{\Pi_1}$ of cubic curves defined by  equation
\begin{equation}\label{node1} z_1z_2z_3+ a_{3,0}z_1^{3}+a_{0,3}z_2^3+a_{0,0}z_3^3=0 \end{equation}
and its projection on the affine space $\text{pr}_1(\mathcal C_{\Pi_1}) =\Pi_1\simeq \mathbb C^3$ in $\mathbb C^9_{1,1}\subset\mathbb P^{9}$.  Note that $(a_{3,0},a_{0,3},a_{0,0})$ are affine coordinates in $\Pi_1$ and everyone can easily check that Hessian
$H_{\mathcal C_{\Pi_1}}$ of this family is given by equation
\begin{equation} \label{node2}
(6^3a_{3,0}a_{0,3}a_{0,0}+2)z_1z_2z_3 +6(a_{3,0}z_1^3+a_{0,3}z_2^3+a_{0,0}z_3^3)=0.\end{equation}
Therefore $h^{-1}(\Pi_1)$ in $\Pi_1\times\mathbb P^2$ is given by
$$z_1z_2z_3=a_{3,0}z_1^3+a_{0,3}z_2^3+a_{0,0}z_3^3=0.$$
Denote $p=(\delta,\delta,\delta)\in \Pi_1$, where $0<\delta\ll 1$, and let $C_p=\kappa^{-1}(p)$. Note that the cubic curve $C_p$ has the same inflection points as Fermat cubic has (see subsection \ref{Hes}).

From the results set forth in subsection \ref{strat}, it follows that locally at the point $\overline a_{3,1}$ with coordinates $a_{3,0}=a_{0,3}=a_{0,0}=0$ the space $\Pi_1$ intersects transversally with irreducible branches of the hypersurface $\mathcal B$ in some neighborhood of this point and the surface $\mathcal B\cap \Pi_1$ is defined in $\Pi_1$ by equation $a_{3,0}a_{0,3}a_{0,0}=0$.
According to Theorem \ref{loc1} the group $\pi_1^{loc} (\mathcal B,\overline a_{3,1})$ is generated by elements
$\gamma_j\in \pi_1(\Pi_1\setminus B,p)$, $j=1,2,3$, represented (see Claim \ref{ex2}) by loops
$$\begin{array}{l} c_1(t)=\{(\delta e^{2\pi it}, \delta,\delta)\in \Pi_1\mid t\in[0,1]\}, \\
c_2(t)=\{(\delta,\delta e^{2\pi it}, \delta)\in \Pi_1\mid t\in[0,1]\}, \\
c_3(t)=\{(\delta,\delta,\delta e^{2\pi it})\in \Pi_1\mid t\in[0,1]\}. \end{array}$$

The preimage $h^{-1}(c_1)$ consists of nine paths
$$\tiny \begin{array}{lll}
q_1(t)=(c_1(t),(0,1,-1)), & q_2(t)=(c_1(t),(1,0,-\lambda(t))), & q_3(t)=(c_1(t),(1,-\lambda(t),0)), \\
q_4(t)=(c_1(t),(0,1,-\omega)), & q_5(t)=(c_1(t),(1,0,-\lambda(t)\omega^2)), & q_6(t)=(c_l(t),(1,-\lambda(t)\omega,0)), \\ q_7(t)=(c_1(t),(0,1,-\omega^2)), &  q_8(t)=(c_1(t),(1,0,-\lambda(t)\omega)), &
    q_9(t)=(c_1(t),(1,-\lambda(t)\omega^2,0)), \end{array} $$
where where $\omega=e^{2\pi i/3}$ and $\lambda(t)=e^{2\pi i t/3}$, and therefore,
$h_*(\gamma_1)=(2,8,5)(3,6,9)=g_2\in \mathcal G_{{3,1}}\subset \mathbb S_9$.

The preimage $h^{-1}(c_2)$ consists of nine paths
$$\tiny \begin{array}{lll}
q_1(t)=(c_2(t),(0,1,-\lambda(t))), & q_2(t)=(c_2(t),(1,0,-1)), & q_3(t)=(c_2(t),(1,-\lambda^2(t),0)), \\
q_4(t)=(c_2(t),(0,1,-\lambda(t)\omega)), & q_5(t)=(c_2(t),(1,0,-\omega^2)), & q_6(t)=(c_2(t),(1,-\lambda^2(t)\omega,0)), \\ q_7(t)=(c_2(t),(0,1,-\lambda(t)\omega^2)), &  q_8(t)=(c_2(t),(1,0,-\omega)), &
    q_9(t)=(c_2(t),(1,-\lambda^2(t)\omega^2,0)) \end{array} $$
and therefore, $h_*(\gamma_2)=(1,4,7)(3,9,6)=g_3$.

Similarly, the preimage $h^{-1}(c_3)$ consists of nine paths
$$\tiny \begin{array}{lll}
q_1(t)=(c_3(t),(0,1,-\lambda^2(t))), & q_2(t)=(c_3(t),(1,0,-\lambda^2(t))), & q_3(t)=(c_3(t),(1,-1,0)), \\
q_4(t)=(c_3(t),(0,1,-\lambda^2(t)\omega)), & q_5(t)=(c_3(t),(1,0,-\lambda^2(t)\omega^2)), & q_6(t)=(c_2(t),(1,-\omega,0)), \\ q_7(t)=(c_3(t),(0,1,-\lambda^2(t)\omega^2)),  & q_8(t)=(c_3(t),(1,0,-\lambda^2(t)\omega)), &
    q_9(t)=(c_2(t),(1,-\omega^2,0)) \end{array} $$
and therefore, $h_*(\gamma_3)=(1,7,4)(2,5,8):=g_4$.
It is easy to see that $g_2g_3=g_4^{-1}$. Therefore the local monodromy group of morphism $h$ at
$\overline a_{3,1}$
is  $\mathcal G_{{3,1}}\simeq\mathbb Z_3^2$ generated in $\mathcal G=Hes$ by permutations $g_2$ and $g_3$.

To prove  statements $(i)$ and $(ii)$, we can assume that $\overline a_{2,1}\in \Pi_1$ has coordinates $a_{3,0}=a_{0,3}=0$ and $a_{0,0}=\varepsilon$, where $0<\epsilon\ll \delta$, i.e., the nodal cubic $C_{\overline a_{2,1}}$ is given by equation
$z_1z_2z_3 + \epsilon z_3^3=0$. Then it is easy to see that the local fundamental group of morphism $h$ at $\overline a_{2,1}$
is generated by elements $\gamma_2$ and $\gamma_3$. And in the case $(i)$ we can assume that the point $\overline a_{1}\in \Pi_1$ has coordinates $a_{3,0}=0$ and $a_{0,0}=a_{0,3}=\varepsilon$, where $0<\epsilon\ll \delta$, i.e. the nodal cubic $C_{\overline a_{2,1}}$ is given by the equation $z_1z_2z_3+\epsilon z_2^3+\epsilon z_3^3=0$. Then it is easy to see that the local fundamental group of morphism $h$ at $\overline a_{1}$ is generated by $\gamma_2$. \qed

\subsection{On the local monodromy groups of the morphism $h$ at points corres\-pond\-ing to cuspidal cubic curves.} \label{casp-loc}
\begin{prop}\label{loc-cusp} The local monodromy groups $\mathcal G_{\overline a}$ of morphism $h$ at points $\overline a\in \mathcal B_{2,2}$ are conjugated to the group $Hes^1$.
\end{prop}
\proof Without loss of generality, we can assume that the cuspidal cubic curve $C_{\overline a_0}\subset \mathbb P^2$ is given by  equation $z_1^3+z_2^2z_3=0$. Its singular point $s$ has coordinates $(0,0,1)$ and the point $p=(0,1,0)$ is its unique nonsingular inflection point.

The fibre $h^{-1} (\overline a_0)$ of morphism $h$ is the intersection of the fibre $C_{\overline a_0}$ of morphism
$\kappa:\mathcal C\to\mathbb P^9$ and Hessian variety $\mathcal H$. Calculating Hessian of the polynomial $z_1^3+z_2^2z_3$, we find that its Hessian curve $H_{C_{\overline a_0}}\subset \mathbb P^2$ given by equation $z_1z_2^2=0$, and it is easy to see that the  intersection numbers of the curves $C_{\overline a_0}$ and $H_{C_{\overline a_0}}$ at points $p$ and $s$ are respectively
$(C_{\overline a_0},H_{C_{\overline a_0}})_p=1$ and $(C_{\overline a_0},H_{C_{\overline a_0}})_s=8$. It follows that in some analytic neighborhood $V\subset \mathcal I$ of the point $(\overline a_0,p)$ the variety $\mathcal I$ is nonsingular and the restriction of the morphism $h$ to $V$ is a biholomorphic isomorphism of the neighborhood $V$ and its image $h(V)$. Therefore, at least one of the orbits of the action of the local monodromy group $\mathcal G_{\overline a_0}$ of morphism $h$ on the set $\{ q_1,\dots, q_9\}$ consists of one point and, hence, $\mathcal G_{\overline a_0}$ is contained in one of the groups conjugated to the group $Hes^1$.

In subsection \ref{gener} it will be shown that the point $\overline a_0\in \mathcal B\cap\Pi$ of the intersection of the variety $\mathcal B$ with a generic projective plane $\Pi\subset \mathbb P^9$ passing through the point $\overline a_0$ is an ordinary cusp of the curve $\mathcal B\cap\Pi$. The group $\pi_1^{loc} (\mathcal B\cap\Pi,\overline a_0)$ (see Claim \ref{Exam-c}) is generated by the geometric generators $\gamma_1$ and $\gamma_2$ satisfying the relation $\gamma_1\gamma_2\gamma_1= \gamma_2\gamma_1\gamma_2$, and the elements
$h_*(\gamma_1)$ and $h_*(\gamma_2)$ generate the group $\mathcal G_{\overline a_0}$.
From the properties of the group $Hes^1$ (see subsection \ref{Hes}), it follows that if $h_*(\gamma_1)\neq h_*(\gamma_2)$, then the group $\mathcal G_{\overline a_0}$ is conjugated in $\mathcal G$ to the group $Hes^1$, or if $h_*(\gamma_1)=h_ * (\gamma_2)$ then $\mathcal G_{\overline a_0}$ is the cyclic group of third order.

Let us  show that the second alternative has no place. To do this, consider the family of  cubic curves $C_{\overline a_{\tau}}$ parameterized by the points of the circle
$S_1=\{ \overline a_{\tau}\in \mathbb P^9 
\}$, $\tau=\delta e^{2\pi i t}$, and given by equation
\begin{equation} \label{kub1} C_{\overline a_{\tau}}: \quad z_1^3+z_2^2z_3+\tau z_3^3=0,\end{equation}
where $t\in [0,1]$ and  $\delta$ is some positive real number. It is easy to see that Hessian variety
$H_{C_{\overline a_{\tau}}}$ of this family is given by equation
\begin{equation}\label{kub2} z_1(3\tau z_3^2- z_2^2)=0.\end{equation}
It follows from (\ref{kub1}) and (\ref{kub2}) that the preimage $h^{-1}(S_1)\subset \mathcal I$ of the circle $S_1$ consists of nine paths three of which are
$$\begin{array}{ll} l_1(t)= & \{ (\overline a(t),\overline z(t))\in \mathcal I\mid \overline a(t)=\overline a_{\tau(t)}, \, \, \overline z(t)=(0,1,0)\} , \\
l_j(t)= & \{ (\overline a(t),\overline z(t))\in \mathcal I\mid \overline a(t)=\overline a_{\tau(t)}, \, \, \overline z(t)=(0,\sqrt[2]{-\delta}e^{\pi i(t+j)},1)\},\quad j=2,3, \end{array}$$
and the other six are
$$l_{j,k}(t)=  \{ (\overline a(t),\overline z(t))\in \mathcal I\mid \overline a(t)=\overline a_{\tau(t)}, \, \, \overline z(t)=(-\sqrt[3]{4\delta}e^{2\pi i(t+j)/3}, \sqrt[2]{3\delta}e^{\pi i(t+k)},1)\},$$
где $j=0,1,2$ и $k=0,1$.
If $\delta \ll 1$, then  $S_1$ (up to conjugation) represents some element $\gamma\in \pi_1^{loc}(\mathcal B,\overline a_0)$ and it is easy to see that the cyclic permutation type of $h_*(\gamma)\in \mathcal G_{\overline a_0}$ is $(6,2,1)$. \qed

\subsection{On local monodromy groups of morphism $h$ at the points
belonging to strata of variety $\mathcal B$ of codimension $\geq 3$.}
Consider the point $\overline a\in \mathcal B_{3,2}$. In any neighborhood of this point there are points $\overline a_1\in \mathcal B_{2,1}$ and $\overline a_2\in \mathcal B_{2,2}$. Therefore, according to Claim \ref{equi-str}, the groups $\mathcal G_{\overline a_1}$ and $\mathcal G_{\overline a_2}\simeq Hes^1$ are subgroups of the local monodromy group $\mathcal G_{\overline a}\subset \mathbb S_9$. The action of the group $\mathcal G_{\overline a_1}$ on the set $\{ 1,\dots, 9\}$ has no fixed elements, and the group $Hes^1$ acts transitively on the set $\{ 2,\dots, 9\}$. Hence, $\mathcal G_{\overline a}$ contains all groups  $Hes^i$, $i=1,\dots,9$.

Since $\mathcal B$ is an irreducible variety, it follows from Proposition \ref{propx} that the monodromy group $\mathcal G$ of the morphism $h$ is generated by permutations of cyclic type $(3,3,1,1,1)$. Each of these generators is contained in some group $Hes^i$. From this it follows
\begin{prop}\label{loc3,2} The local monodromy groups $\mathcal G_{\overline a}$ of the morphism $h$ at the points $\overline a\in\mathcal B_{3,2}$ and the local monodromy groups $\mathcal G_{\overline a}$ at the points $\overline a$, belonging to the strata of $\mathcal B'$ of the variety $\mathcal B$ of codimension $\geq 4$, coincide with the whole group $\mathcal G=Hes$.
\end{prop}

\section{The surfaces of inflection points of generic two-dimensional linear systems of plane cubic curves}\label{two}
\subsection{Two-dimensional linear systems of plane cubic curves.} Let $\overline z=(z_1,z_2,z_3)$ be homogeneous coordinates in $\mathbb P^2$ and $\overline t=(t_1,t_2,t_3)$ homogeneous coordinates in a projective plane $\Pi\subset \mathbb P^9$. Consider (projectively) two-dimensional linear system of plane cubic curves $C_{\Pi}$ given by equation $F(\overline t,\overline z)=0$, where
$$F(\overline t,\overline z)=t_1F(\overline a_1,\overline z)+t_2F(\overline a_2,\overline z)+t_3F(\overline a_3,\overline z).$$
Let $H(\overline t,\overline z)= \det (\frac{\partial^2F(\overline t,\overline z)}{\partial z_i\partial z_j})$ be the Hessian of the polynomial $F(\overline t,\overline z)$. The equations
\begin{equation}\label{equations} F(\overline t,\overline z)=0, \quad H(\overline t,\overline z)=0\end{equation}
define in $\Pi\times \mathbb P^2$ a surface  $\mathcal I_{\Pi}=\mathcal I\cap \text{pr}^{-1}(\Pi)$ of the inflection points of the linear system of cubic curves $C_{\Pi}$.

Set $x=\frac{z_1}{z_3}$, $y=\frac{z_2}{z_3}$ and $\alpha=\frac{t_2}{t_1}$, $\beta=\frac{t_3}{t_1}$. Below, each time, without specifying, we will choose the base elements $F_i(\overline z):=F(\overline a_i,\overline z)$ of the linear system $\mathcal C_{\Pi}$ and homogeneous coordinates
$(z_1,z_2,z_3)$ in $\mathbb P^2$ so that the form of equations (\ref{equations}) in non-homogeneous coordinates $(\alpha, \beta,x,y)$ in $\mathbb C^2\times \mathbb C^2\subset \Pi\times\mathbb P^2$ would be most convenient to study the properties of the variety
$\mathcal I_{\Pi}$.

\subsection{Generic projections to the plane.} \label{gener}
The complete linear system of cubic curves, parametrized by the points of the space $\mathbb P^9$, defines, so called, Veronese-3 imbedding $\varphi_3:\mathbb P^2\hookrightarrow \hat{\mathbb P}^9$, $\varphi_3(\mathbb P^2)=\mathcal V_3$, of the projective plane $\mathbb P^2$ into the dual projective space of the parameter space $\mathbb P^9$. Two-dimensional linear subspaces $\Pi\subset \mathbb P^9$ define linear projections $\text{pr}_{\Pi}:\hat{\mathbb P^9}\to\hat{\Pi}$ into the projective planes $\hat{\Pi}$ dual to the projective planes $\Pi$. It is well known (see, for example, \cite{C-F}) that for planes  $\Pi$ belonging to a Zariski open everywhere dense set $\mathcal W$ in Grassmannian $Gr(3,10)$ of two-dimensional planes of $\mathbb P^9$, the compositions
 $\xi :=\text{pr}_{\Pi}\circ\varphi_3:\mathbb P^2\to \hat{\Pi}$ are generic coverings of the projective plane, i.e., they have the following properties:
\begin{itemize}
\item[$(i)$] $\xi$ is a finite morphism of degree nine;
\item[$(ii)$] the ramification curve $\widetilde R \subset \mathbb P^2$ of $\xi$ is non-singular and ramification index along $\widetilde R$ is equal to two;
\item[$(iii)$] the discriminant  $\hat B=\xi(\widetilde R)$ has only ordinary nodes and cusps as its singular points and the restriction of $\xi$ to $\widetilde R$ is a birational morphism.
\end{itemize}

It follows from the proof of Theorem 4 in \cite{K-C} that the degree of the curve $\hat B$ is equal to $18$, and the curve
$\hat B$ has $42$ ordinary cusps, $84$ ordinary nodes, and the geometric genus $g(\hat B)$ of
$\hat B$ is equal to $10$. The curve $\widetilde R\subset \mathbb P^2$ is a smooth plane curve of degree six.

The cubic curves $C_{\overline a}=\xi^{-1}(L_{\overline a})$, parametrized by the points of $\Pi$, are inverse images of lines $L_{\overline a}$ in $\hat{\Pi}$ dual
to the points $\overline a\in \Pi$. A cubic curve $C_{\overline a}$ is singular if and only if  $L_{\overline a}$ is a tangent line of the curve
$\hat B$. Therefore the curve $\hat B$ is dual to the curve $B:=\mathcal B\cap \Pi$ and if a cubic curve $C_{\overline a}$ is singular, then its singular points belong to $C_{\overline a}\cap \widetilde R$.  For a point $q\in L_{\overline a}\cap \hat B$, the only one point from $\xi^{-1}(q)$ can be a singular point of the curve $C_{\overline a}$. In addition, the restriction of projection $\text{pr}_2: \Pi\times\mathbb P^2\to\mathbb P^2$ to $S=\mathcal S\cap (\Pi\times \mathbb P^2)$ is a one-to-one map onto the curve $\widetilde R$.

It follows from equisingular stratification of singular cubic curves (see subsection \ref{strat}) that $\hat B$ can have only bitangent and three-tangent lines   $L_{\overline a_i}$ (in this case cubic curves $C_{\overline a_i}$ are unions of transversally intersecting lines and quadrics, i.e.,
$\overline a_i\in \mathcal B_{2,1}$, and, respectively, they are unions of three lines, i.e., $\overline a_i\in \mathcal B_{3,1}$), and if $L_{\overline a_i}$ is the tangent line of the curve $\hat B$ at its inflection point, then $C_{\overline a_i}$ is a cuspidal cubic, i.e., $\overline a_i\in \mathcal B_{2,2}$. Therefore, by Pl$\ddot{u}$kker formulae (see, for example, \cite{Ku}), the degree of the curve $B$ is equal to $12$ and $B$ has $24$ ordinary cusps, $n_1$ ordinary nodes, and $n_2$ ordinary triple poins (the singularities of type $D_4$) as its singular points, where $n_1+3n_2=21$.
\begin{rem} Note that if the plane $\Pi$ is sufficiently general, i.e. $\Pi\cap \mathcal B_{3,1}=\emptyset$, then the curves $B\subset \Pi$ and $\hat{B}\subset \hat{\Pi}$ give an example of two dual cuspidal curves for which the fundamental groups  $\pi_1(\Pi\setminus B)$ and $\pi_1(\hat{\Pi}\setminus \hat B)$ are not abelian groups.
\end{rem}

Denote
$$\mathfrak{B}_{\widetilde R,1}=\{ \overline z\in \widetilde R\mid \xi(\overline z)\, \,
\text{is an inflection point of}\,\, \hat B\},$$
$$\mathfrak{B}_{\widetilde R,2}=\{ \overline z\in \widetilde R\mid \xi(\overline z)\, \, \text{is a cusp of}\,\, \hat B\}.$$
The set  $\mathfrak{B}_{\widetilde R,1}$ consists of $24$ points and $\mathfrak{B}_{\widetilde R,2}$ consists of $42$ points.

Further we assume that $\Pi \in \mathcal W$.
\subsection{Properties of Stein coverings associated with  morphisms $h_{\Pi}$.}\label{Stein}
A two-dimensional linear system of cubic curves $\mathcal C_{\Pi}$ defines a generic covering of the plane $\xi:\mathbb P^2\to \hat{\Pi}$ and it follows from van Kampen -- Zariski Theorem that the monodromy group of the covering $h_{\Pi}:\mathcal I_{\Pi}\to \Pi$ is the group $\mathcal G=Hes\subset \mathbb S_9$.

Let us consider the Stein covering $\varphi :Y\to \Pi$,  $\deg \varphi=9$ associated with $h_{\Pi}$. The covering $\varphi$ is unramified over $\Pi\setminus B$. The following description of properties of the covering $\varphi: \varphi^{-1}(V)\to V$, where $V\subset \Pi$ is a sufficiently small complex-analytic neighbourhood of a point  $o\in B$, based on results of section \ref{mon-loc} and subsection \ref{monodr}.

If $o\in \mathcal B_1$, then the local monodromy group of the covering $\varphi$ at $o$ is the cyclic group generated by permutation whose cyclic type is $(1,1,1,3,3)$. Therefore $\varphi^{-1}(V)$ is a disjoint union of five irreducible open in $Y$ sets, $\varphi^{-1}(V)=\cup_{i=1}^{5}U_i$, each of which is non-singular and, up to numbering, the coverings $\varphi: U_i\to V$, $i=1,2,3$, are bi-holomorphic maps and  $\varphi:U_i\to V$, $i=4,5$, are three-sheeted coverings branched along $B\cap V$.

If $o\in \mathcal B_{2,1}$ and $B\cap V= B'_1\cup B'_2$, where $B'_1$ and $B'_2$ are two branches of the curve $B$ intersecting transversally at $o$, then the local monodromy group $\mathcal G_{2,1}\subset \mathbb S_9$ of the covering $\varphi$ at $o$ is the abelian group generated by permutations $g_3=(2,8,5)(3,6,9)$ (the image of bypass around $B'_1$) and $g_4=(1,4,7)(3,9,6)$ (the image of bypass around $B'_2$). The set $\{ 1,\dots, 9\}$ splits into three orbits  $\{ 1,4,7\}$, $\{ 2,8,5\}$, and $\{ 3,9,7\}$ of the action of the group $\mathcal G_{2,1}$. Therefore $\varphi^{-1}(V)$ is a disjoint union of three irreducible open in $Y$ sets, $\varphi^{-1}(V)=\cup_{i=1}^{3}U_i$.  Up to numbering, the sets $U_1$ and $U_2$ are non-singular and $U_3$ has a singularity at the point
$o'= U_3\cap \varphi^{-1}(o)$. The coverings $\varphi: U_i\to V$, $i=1,2,3$, are three-sheeted, the coverings $\varphi:U_1\to V$  is ramified over $B'_{1}$, the covering $\varphi:U_2\to V$ is ramified over $B'_{2}$, and the covering $\varphi:U_3\to V$ is ramified over $B'_{1}\cup B'_2$. Note (see the proof of Claim \ref{deg3}) that if $\psi: W\to U_3$ is the minimal desingularisation of the singular point $o'$, $E=\psi^{-1}(o')$, then $E$ is the smooth rational curve and $(E^2)_W=-3$.

If $o\in \mathcal B_{3,1}$ and $B\cap V= B'_1\cup B'_2\cup B'_3$, where $B'_1$, $B'_2$, and $B'_3$ are three branches of the curve $B$ pairwise meeting transversally at $o$, then the local monodromy group $\mathcal G_{3,1}\subset \mathbb S_9$ of the covering $\varphi$ at $o$ is the abelian group generated by three permutations $g_3=(2,8,5)(3,6,9)$ (the image of bypass around $B'_1$),  $g_4=(1,4,7)(3,9,6)$ (the image of bypass around $B'_2$), and $g_3=(2,5,8)(1,7,4)$ (the image of bypass around $B'_3$). The set  $\{ 1,\dots, 9\}$ splits into three orbits  $\{ 1,4,7\}$, $\{ 2,8,5\}$, and $\{ 3,9,7\}$ of the action of the group $\mathcal G_{3,1}$.  Therefore $\varphi^{-1}(V)$ is a disjoint union of three irreducible open in $Y$ sets, $\varphi^{-1}(V)=\cup_{i=1}^{3}U_i$.  Up to numbering, the sets $U_1$, $U_2$, and $U_3$ have a singularity at the points $o'_i= U_i\cap \varphi^{-1}(o)$. The coverings $\varphi: U_i\to V$, $i=1,2,3$,  are three-sheeted, the covering $\varphi:U_1\to V$ is ramified over $B'_{1}\cup B'_2$, the covering $\varphi:U_2\to V$ is ramified over $B'_{2}\cup B'_3$, and the covering $\varphi:U_3\to V$ is ramified over $B'_{1}\cup B'_3$. As above, note that if $\psi: W_i\to U_i$ are the minimal desingularisations of the singular points $o'_i$, $E_i=\psi^{-1}(o'_i)$, then $E_i$ are smooth rational curves and $(E_i^2)_{W_i}=-3$.

If $o\in \mathcal B_{2,2}$, then the group of local monodromy $\mathcal G_{2,2}\subset \mathbb S_9$ of the covering $\varphi$ at $o$ is the group
$\mathcal G_{2,2}\simeq Hes^1$ generated by two permutations $g_1=(4,5,6)(7,9,8)$ and $g_2=(2,8,5)(3,6,9)$. The set $\{ 1,\dots, 9\}$ splits into two orbits   $\{ 1\}$ and $\{ 2,\dots, 9\}$ of the action of the group $\mathcal G_{2,2}$.  Therefore $\varphi^{-1}(V)$ is a disjoint union of two irreducible open in $Y$ sets, $\varphi^{-1}(V)=U_1\cup U_2$. The covering
$\varphi: U_1\to V$ is a bi-holomorphic map and therefore $U_1$ is non-singular. The degree of covering $\varphi : U_2\to V$ is equal to eight and this covering is branched over the curve   $B\cap V= B'$ having the singularity of type $A_2$ at $o$. The group $Hes^1\subset \mathbb S_8$ is the monodromy group of the covering $\varphi: U_2\to V$ ($\mathbb S_8$ acts on $\{ 2,\dots, 9\}$).

Let us show that $U_2$ is nonsingular, the ramification curve  $R'\subset U_2$ of the covering $\varphi$ is also nonsingular, and the map $\varphi: R'\to B'$ is two-sheeted and it is ramified at the point $o$. For this, let us remind that the group
$Hes^1=\langle g_1,g_2\rangle$ can be imbedded into $GL(2,\mathbb C)$ (the group no. 4 in \cite{S-T}, see also \cite{Ku3}) and under this imbedding the elements $g_1$ and $g_2$ are  3-reflections:
$$ g_1=\frac{-\omega}{\sqrt{2}}\left(\begin{array}{cc} \varepsilon & \varepsilon^3 \\ \varepsilon & \varepsilon^7\end{array}\right), \quad
g_2=\frac{\omega}{\sqrt{2}}\left(\begin{array}{cc} \varepsilon^3 & \varepsilon^5 \\ \varepsilon^7 & \varepsilon^5\end{array}\right), $$
where $\omega=e^{2\pi i/3}$ and $\varepsilon= e^{2\pi i/8}$. Since the group $Hes^1\subset GL(2,\mathbb C)$ is generated by reflections, the factor-space  $\mathbb C^2/Hes^1$ is  $\mathbb C^2$. Denote by $\widetilde{\psi}:\mathbb C^2\to\mathbb C^2$ the factorization morphism.
We have $\deg \widetilde{\psi}=24$ and, as it is known (see \cite{Ku3}),  $\widetilde{\psi}$ is defined by functions $\Psi=u^4-2\sqrt{-3}u^2v^2+v^4$,
$\Theta=uv(u^4-v^4)$ and it is ramified with multiplicity three along four lines $\cup_{i=1}^4L_i$ given by equation $\Phi=0$, where
$\Phi=u^4+2\sqrt{-3}u^2v^2+v^4$. The branch curve $\widetilde B$ is given by $\Psi^3+12\sqrt{-3}\Theta^2=0$.
The degree of the restriction of $\widetilde{\psi}$ to each line $L_i$, $\widetilde{\psi}: L_i\to \widetilde B$, equals two. The group $\langle g_1\rangle$ transpose the lines $L_2$, $L_3$, and $L_4$ and leaves fixed the line $L_1$ given by equation $(\sqrt{2}+\omega\varepsilon)u+\omega\varepsilon^3v=0$. The covering $\widetilde{\psi}$ is the composition $\widetilde{\psi}=\psi\circ\overline{\psi}$, where
$\overline{\psi}:\mathbb C^2\to\mathbb C^2/\langle g_1\rangle\simeq \mathbb C^2$ is the three-sheeted factorization map defined by the action of the group $\langle g_1\rangle$. It is ramified along $L_1$ and $\psi:\mathbb C^2/\langle g_1\rangle\to\mathbb C^2$ is a eight-sheeted covering ramified with multiplicity three along nonsingular curve $\overline R=\overline{\psi}(L_2\cup L_3\cup L_4)$ and whose monodromy group is $Hes^1\subset \mathbb S_8$ acting on eight left cosets of the group $\langle g_1\rangle$.
The covering $\psi:\overline R\to\widetilde B$ is a two-sheeted map ramified at the origin. Now, to show that $U_2$ is nonsingular, the ramification curve
$R'\subset U_2$ of $\varphi$ is also nonsingular, and  $\varphi: R'\to B'$ is a two-sheeted covering branched at $o$, it suffices to imbed the pair $(V,B')$ into $(\mathbb C^2,\widetilde B)$ and apply Grauert - Remmert - Riemann - Stein Theorem, Proposition \ref{gal-st}, and Claim \ref{norm-cusp} to the coverings
$\psi: \psi^{-1}(V)\to V$ and $\varphi: U_2\to V$.

Denote by $\mathcal R\subset Y$ the ramification curve of $\varphi$. It follows from above that $\varphi:\mathcal R\to B$ is a two-sheeted covering branched at  $24$ points (the cusps of $B$). Therefore, by Hurwitz formula, the geometric genus $g(\mathcal R)$ of $\mathcal R$ is equal to
$$ g(\mathcal R)=\frac{1}{2}[2(2g(B)-2)+24]+1=31.$$

The number of singular points of the surface $Y$ equals $n_1+3n_2=21$.
\begin{prop}\label{St-inv} Let $\sigma:X\to Y$ be the minimal resolutions of singularities of $Y$. The square of the canonical class  $K_X$ of $X$ equals  $K_X^2=18$ and the Euler characteristic  $\chi(\mathcal O_X)=\sum_{i=0}^2(-1)^i\dim H^i(X,\mathcal O_X)$ of the structure sheaf $\mathcal O_X$ of $X$ equals $\chi(\mathcal O_X)=1-q+p_g=9$.
\end{prop}
\proof Denote $\psi=\varphi\circ\sigma$ and denote by $R\subset X$ the proper preimage of the curve
$\mathcal R$. By Claim \ref{deg3},  $\sigma^{-1}(Sing\, Y)$ is the disjoint union of $21$ rational curves,  $E:=\sigma^{-1}(Sing\, Y)=\sqcup_{i=1}^{21}E_i$. In addition, by the same Claim, we have  $(E_i^2)_X=-3$ and $(R,E_i)_{X}=2$.

We have $K_X=\psi^*(K_{\Pi})+2R+E$, since $\psi$ is ramified along $R$ with multiplicity three, $(E_i,\psi^*(K_{\Pi})_X=0$ and $(K_X+E_i,E_i)_X=-2$.
In addition, we have
$(\psi^*(K_{\Pi}),\psi^*(K_{\Pi}))_X=81$, since $\deg \psi=K^2_{\Pi}=9$ , and $(R,\psi^*(K_{\Pi}))_X =-72$, since the degree of the covering
$\psi:R\to B$ equals two and $\deg B=12$.
Therefore the  equalities
\begin{equation}\label{R-sq}(R^2)_X= (R,K_X)_X=30\end{equation}
follow from equalities  $(K_X+R,R)_X= (\psi^*(K_{\Pi})+3R+E,R)_X=60$. As a result, we have
$K_X^2=(\psi^*(K_{\Pi})+2R+E,\psi^*(K_{\Pi})+2R+E)_X=18$.

Denote by $e(M)$ Euler topological characteristic of a space $M$. We have
$e(X)=9e(\Pi\setminus B)+5e(B\setminus Sing\, B)+4n_1+6n_2+2c$,
since $\deg \psi =9$, the number of preimages of each nonsingular point of the curve $B$ equals five, the preimage of each node of  $B$ is a disjoint union of projective line ($e(\mathbb P^1)=2$) and two points, the preimage of each triple point of $B$ is a disjoint union of three projective lines, and the preimage of each cusp of $B$ consists of two points. And since $g(B)=10$, $n_1+3n_2=21$, $c=24$, and $e(\Pi)=3$, then
$e(B)=-2(g(B)-1)-n_1-2n_2=n_2-39$, $e(\Pi\setminus B)=42-n_2$, and
$e(B\setminus Sing\, B)=-2(g(B)-1)-2n_1-3n_2-c=3n_2-84$. Therefore, $e(X)=90$ and applying Noether formula, we obtain that $\chi(\mathcal O_X)=\frac{K_X^2+e(X)}{12}=9$. \qed

\subsection{Resolution of singularities of a surface of inflection points of generic two-dimensional linear system of plane cubic curves.}\label{sing}
Denote by  $\nu: I:= I_{\Pi}\to \mathcal I_{\Pi}$ the normalization of variety $\mathcal I_{\Pi}$ and let  $p_1=h_{\Pi}\circ \nu: I\to\Pi$. Denote also  $S=\mathcal S\cap \mathcal C_{\Pi}$. By Claim \ref{sing-inf}, the curve $S$ is contained in $\mathcal I_{\Pi}$.
\begin{prop}\label{main1} The surface $I$ is nonsingular and the morphism $p_1: I\to\Pi$ coincides with $\psi=\varphi\circ\sigma: X\to \Pi$, where
$\varphi :Y\to \Pi$ is the Stein covering associated with  $h_{\Pi}$ and  $\sigma: X\to Y$ is the minimal resolution of singularities of the surface $Y$.
\end{prop}
\proof
Obviously,  $\mathcal  I_{\Pi}$ is nonsingular at the points belonging to $h_{\Pi}^{-1}(\Pi\setminus B)$.
\begin{claim}\label{non-sing} The surface $\mathcal I_{\Pi}$ is nonsingular at the points belonging to
$\mathcal I_{\Pi}\setminus  S$.
\end{claim}
\proof An irreducible nodal cubic $C_{\overline a_1}$, $\overline a_1\in \mathcal B_1\cap \Pi$, has three smooth inflection points (respectively, a cuspidal cubic $C_{\overline a_1}$,
$\overline a_1\in \mathcal B_{2,2}\cap \Pi$, has a single smooth inflection point). It was shown in \cite{Ku} that the intersection number of the cubic curve
$C_{\overline a_1}$ and its Hesse curve $H_{\overline a_1}$ at its singular point $s\in C_{\overline a_1}$ equals $(C_{\overline a_1},H_{\overline a_1})_s=6$ (respectively, $(C_{\overline a_1},H_{\overline a_1})_s=8$).
Therefore, by Claim \ref{connect}, the surface $\mathcal I_{\Pi}$ is nonsingular at the points belonging to $\mathcal I_{\Pi}\setminus (S\cup h^{-1}(\mathcal B_{2,1}\cup\mathcal B_{3,1})$. Therefore, to complete the proof of Proposition \ref{non-sing}, it suffices to show that  $\mathcal I_{\Pi}$ is non-singular at the points of $h^{-1}(\overline a_1)\setminus \mathcal S$ for $\overline a_1\in (\mathcal B_{2,1}\cup\mathcal B_{3,1})\cap \Pi$.

 Without loss of generality, we can assume that $z_1z_2z_3+a_{3,0}z_1^3=0$ is an equation of the cubic curve $C_{\overline a_1}$ (if
 $\overline a_1\in \mathcal B_{3,1}$, then $a_{3,0}=0$ and $h_{\Pi}^{-1}(\overline a_1)$ are three lines in $\{ \overline a_1\}\times\mathbb P^2$ given by  equations $z_i=0$, $i=1,2,3$, and if $\overline a_1\in \mathcal B_{2,1}$, then $a_{3,0}\neq 0$ and $h_{\Pi}^{-1}(\overline a_1)$ is the line in
 $\{ \overline a_1\}\times\mathbb P^2$ given by equation $z_1=0$).

Let us show that the points of line $E=\{ \overline a_1\}\times\{ z_1=0\}$ are non-singular points of  $\mathcal I_{\Pi}$ if $z_2z_3\neq 0$ (in the case when  $\overline a_1\in \mathcal B_{3,1}$  and the points belong to the lines $\{ z_2=0\}$ and $\{ z_3=0\}$ not being singular points of the curve
$C_{\overline a_1}$, the proof that the surface $\mathcal I_{\Pi}$ is non-singular at these points is similar and therefore it will be omitted). For this, let us choose base elements $F_1(\overline z)$, $F_2(\overline z)$, $F_3(\overline z)$ of linear system $\mathcal C_{\Pi}$ so that in non-homogeneous coordinates  $(\alpha,\beta,x,y)$  the linear system $\mathcal C_{\Pi}$ is given by equation
\begin{equation}\label{eq01} (xy+a_{3,0}x^3)+\alpha(\sum_{1\leq i+j\leq 3}b_{i,j}x^iy^j)+
\beta(1+\sum_{1\leq i+j\leq 3}c_{i,j}x^iy^j)=0.\end{equation}
Denote by $\mathfrak{m}$  the ideal in the polynomial ring $\mathbb C[\alpha,\beta,x,y]$  generated by monomials  $\alpha$, $\beta$, and $x$ and write equation (\ref{eq01}) in the following form:
\begin{equation}\label{eq02} xy+\alpha(\sum_{1\leq j\leq 3}b_{0,j}y^j)+
\beta(1+\sum_{1\leq j\leq 3}c_{0,j}y^j) +R_1(\alpha,\beta,x)=0,\end{equation}
where $R_1(\alpha,\beta,x)\in \mathfrak{m}^2$. Simple calculations (which we omit) show that in non-homogeneous coordinates $(\alpha,\beta,x,y)$ the Hessian of the linear system $\mathcal C_{\Pi}$ has the following form:
\begin{equation}\label{eq03} 2[xy+\alpha(\sum_{1\leq j\leq 2}b_{0,j}y^j-3b_{0,3}y^3)+
\beta(-3+\sum_{1\leq j\leq 2}c_{0,j}y^j -3c_{0,3}y^3)] +R_2(\alpha,\beta,x)=0,\end{equation}
where $R_2(\alpha,\beta,x)\in \mathfrak{m}^2$.

The surface $\mathcal I_{\Pi}\cap (\mathbb C^2\times\mathbb C^2)$ in $\mathbb C^2\times\mathbb C^2$ is given by equations (\ref{eq02}) and (\ref{eq03}), and the line $E$ is given by equations $\alpha=\beta=x=0$. Therefore necessary and sufficient condition for a point  $p=(0,0,0,y_0)$ with $y_0\neq 0$ to be singular point of $\mathcal I_{\Pi}$ is the proportionality of (linear in variables  $\alpha$, $\beta$, and $x$) forms entering into left parts of equations  (\ref{eq02}) and (\ref{eq03}), and in which $y$ takes the concrete value, namely, $y_0$.  It is easy to see that one of necessary conditions is vanishing the coefficient  $b_{0,3}$. In this case the cubic curve $C_{\overline a_2}$ given by equation $F_2(\overline z)=0$ contains not only the singular point $q_1=(0,0,1)$ of the curve $C_{\overline a_1}$, but also it contains the second singular point $q_2=(0,1,0)$ of $C_{\overline a_1}$. In our case it is impossible, since the linear system $\mathcal C_{\Pi}$ defines the generic covering $\xi:\mathbb P^2\to \hat{\Pi}$ of the projective plane. \qed

\begin{claim} \label{non-singS}
The curve $S$ is non-singular at
$s\in S\subset \mathcal I_{\Pi}$ if $h_{\Pi}(s)\in \mathcal B_1\cup \mathcal B_{2,1}\cup \mathcal B_{3,1}$.
\end{claim}
\proof As above, we can assume that the point $s$ lies in neighbourhood
$\mathbb C^2\times \mathbb C^2\subset \Pi\times \mathbb P^2$ with coordinates $(\alpha,\beta,x,y)$ and has coordinates $(0,0,0,0)$. Choose the base polynomials $F_1(x,y)$, $F_2(x,y)$, and $F_3(x,y)$ of linear system $\mathcal C_{\Pi}$ as follows:
$$F_1(x,y)=xy+a_{3,0}x^3+a_{0,3}y^3,\,\, F_2(x,y)=\sum_{1\leq i+j\leq 3}b_{i,j}x^iy^j,\,\,F_3(x,y)=1+\sum_{1\leq i+j\leq 3}c_{i,j}x^iy^j,$$
i.e., the linear system $\mathcal C_{\Pi}$ is given by
\begin{equation}\label{eq04} (xy+a_{3,0}x^3+a_{0,3}y^3)+\alpha(\sum_{1\leq i+j\leq 3}b_{i,j}x^iy^j)+\beta(1+\sum_{1\leq i+j\leq 3}c_{i,j}x^iy^j)=0.\end{equation}
We can assume that $b_{1,1}=c_{1,1}=0$ (correcting  $F_2(x,y)$ and $F_3(x,y)$ with help of the polynomial $F_1(x,y)$ if it is necessary) and $b_{1,0}=1$, $c_{1,0}=0$.
Denote by  $\mathfrak{m}$  the ideal in the power series ring $\mathbb C[[\alpha,\beta,x,y]]$ generated by monomials $\alpha$, $\beta$, $x$, and $y$, and denote by $\mathfrak{m}_{\hat{\beta}}$  the ideal in the power series ring $\mathbb C[[\alpha,x,y]]$ generated by monomials $\alpha$, $x$, and $y$.
Denote also by $w$ the left side of equation (\ref{eq04}) and let $w_1=\frac{\partial w}{\partial x}$ and $w_2=\frac{\partial w}{\partial y}$. We have $w=\beta+R_0(\alpha,\beta,x,y)$ and
\begin{equation}\label{eq08}   w_1=y+\alpha +R_1(\alpha,\beta,x,y),\,\,  w_2=x+b_{0,1}\alpha+ c_{0,1}\beta+R_2(\alpha,\beta,x,y),
\end{equation}
where $R_i(\alpha,x,y)\in \mathfrak{m}^2$ for $i=0,1,2$. It is easy to see from (\ref{eq08}) that the functions $w$, $w_1$, $w_2$, and $\alpha$  are local parameters in some neighbourhood $V$ of the point $s$, the variety $U=\mathcal C_{\Pi}\cap V$ is given in $V$ by equation $w=0$, and the curve $S\cap V$ is given in $U$ by equations $w_1=w_2=0$. \qed

\begin{claim} \label{sin-node} If $h_{\Pi}(s)\in \mathcal B_1\cup \mathcal B_{2,1}\cup \mathcal B_{3,1}$ for
$s\in S\subset \mathcal I_{\Pi}$, then in some analytic neighbourhood of the point $s$ the surface $\mathcal I_{\Pi}$ is the union of two non-singular surfaces meeting transversally at $s$.
\end{claim}
\proof We use notations introduced in the proof of Claim \ref{non-singS}.
Simple calculations (which we omit) show that in non-homogeneous coordinates $(\alpha,\beta,x,y)$ the equation of Hesse surface $H_{\mathcal C_{\Pi}}$ has the following form
\begin{equation}\label{eq05} 2[-3\beta+xy+\alpha(x+b_{0,1}y)+\beta(-3+c_{0,1}y)+4b_{01}\alpha^2+4c_{0,1}\alpha\beta]+R_3(\alpha,\beta,x,y) =0,\end{equation}
where $R_3(\alpha,\beta,x,y)\in \mathfrak{m}^3$.

It follows from (\ref{eq04}) that the restriction of function $\beta$ to $V$ has the following form:
\begin{equation}\label{eq06} \beta= -[xy+\alpha (x+ b_{0,1} y)]+ R_4(\alpha,x,y),\end{equation}
where $R_4(\alpha,x,y)\in \mathfrak{m}^3_{\hat{\beta}}$.
Therefore the surface $\mathcal I_{\Pi}\cap U$ is given by equation (\ref{eq05}) in which the right side of  equation  (\ref{eq06}) is substituted instead of   $\beta$. After this substitution and simple transformations, equation  (\ref{eq05}) has the following form:
\begin{equation}\label{eq07} 4(y+\alpha)(x+b_{0,1}\alpha)+R_5(\alpha,x,y)=0,\end{equation}
where $R_5(\alpha,x,y)\in \mathfrak{m}^3_{\hat{\beta}}$. Therefore, $s$ is a singular point of  $\mathcal I_{\Pi}$ of multiplicity two. Also, it follows from this that $S=Sing\, \mathcal I_{\Pi}$, since the point $s$ is a generic point belonging to $S$.

It follows from (\ref{eq07}) that $s$, as a singular point of the curve $D=\{ \alpha=0\} \cap \mathcal I_{\Pi}$, is an ordinary node.
Let $\sigma :\widetilde U\to U$ be the monoidal transformation with center in the curve $S\cap U$, $\sigma^{-1}(S\cap U)=E$. It is easy to see that the proper preimage $\sigma^{-1}(D)$ of $D$ intersects transversally with $E$ at two points, and since  $s$ is a generic point of the curve $S\cap U$, then the intersection $\sigma^{-1}(U\cap\mathcal I_{\Pi})\cap E$ of the proper preimage of the surface $U\cap\mathcal I_{\Pi}$ with $E$ is a disjoint union of two nonsingular curves if the neighbourhood  $V$ is chosen small enough. It follows from this that the surface $\mathcal I_{\Pi}\cap U$ is the union of two nonsingular surfaces meeting transversally, since the multiplicity of the singularity of the surface $\mathcal I_{\Pi}\cap U$ at each point $s\in S\cap U$ equals two. \qed

\begin{claim} \label{sin-cusp} If $h_{\Pi}(s)\in \mathcal B_{2,2}$ for $s\in S\subset \mathcal I_{\Pi}$, then the surface $I$ is nonsingular at the points of  $\nu^{-1}(s)$.
\end{claim}
\proof Let us consider a sufficiently small complex-analytic neighbourhood $V\subset \Pi$ of the point $p_1(s)$. The morphism  $p_1:p_1^{-1}(V)\to V$ is a finite covering and  $p_1^{-1}(V)$ is a normal surface (consisting of two components). Therefore Claim \ref{sin-cusp} follows from   Grauert - Remmert - Riemann - Stein Theorem and results of subsection \ref{Stein}. \qed \\

To complete the proof of Proposition \ref{main1}, it suffices to apply Claims \ref{non-sing} -- \ref{sin-cusp}, \ref{deg3} and  Grauert - Remmert - Riemann - Stein Theorem. \qed

\section{Calculation of the irregularity}
\subsection{The additional condition of generality.}
Further we use notations introduced in section \ref{two}.
To compute the irregularity $q=\dim H^{1}(\mathcal O_I,I)$ of a surface $I$, 
it is useful to introduce an {\it additional generality condition} of two-dimensional linear systems $\mathcal C_{\Pi}$ of plane cubic curves. We say that a linear system $\mathcal C_{\Pi}$ of cubic curves, parameterized by points of the projective plane $\Pi\in \mathcal W\subset Gr(3,10)$, satisfies the {\it additional condition of generality} if for any point $\overline z\in \mathbb P^2$ there is only a finite number of cubic curves of $\mathcal C_{\Pi}$, passing through $\overline z$ for which $\overline z$ is their inflection point.

Denote by $J_{1, o}$ a subset of planes $\Pi$ in $\mathcal W\subset Gr(3,10)$ such that the point
$o\in \mathbb P^2\setminus \widetilde R$ is an inflection point of each cubic from $\mathcal C_{\Pi}$ containing the point $o$.
\begin{claim} \label{inf1} We have $\dim J_{1,o}\leq 17$.
\end{claim}
\proof Since $o\in\mathbb P^2\setminus \widetilde R$, then in some neighborhood of $o$ the generic covering $\xi:\mathbb P^2\to \hat{\Pi}$ is a biholomorphic isomorphism. Therefore
the base elements $F_1 (x,y)$, $F_2 (x,y)$, and $F_3 (x,y)$ of the linear system $\mathcal C_{\Pi}$ of the following form:
$$F_1(x,y)=x+\sum_{2 \leq i+j\leq 3}b_{i,j}x^iy^j, \,\, F_2(x,y)= y+
\sum_{2 \leq i+j\leq 3}c_{i,j}x^iy^j,$$ $$F_3(x,y)=1+\sum_{2 \leq i+j\leq 3}d_{i,j}x^iy^j$$
are defined uniquely by the point $o\in \mathbb P^2\setminus \widetilde R$ for any plane
$\Pi\in \mathcal W$. Obviously, the necessary and sufficient condition for $\Pi$ to belong to $J_{1,o}$ is the divisibility of the quadratic form
$$ (t_1b_{2,0}+t_2c_{2,0})x^2 +(t_1b_{1,1}+t_2c_{1,1})xy+(t_1b_{0,2}+t_2c_{0,2})y^2 $$
by linear form $t_1x + t_2y$ for each point $(t_1,t_2)\in \mathbb P^1$.
It is easy to see (putting $x=t_2=1$) that this condition is equivalent to the following equalities:
$$b_{2,0}=c_{1,1},\quad b_{1,1}=c_{0,2},\quad b_{0,2}=c_{2,0}=0.$$
Consequently, the dimension of the variety $J_{1, o}$ does not exceed the total (equals to $17$) number of undefined coefficients $b_{I,j}$, $c_{I, j}$ and $d_{I,j}$ for monomials in polynomials  $F_1(x,y)$, $F_2(x,y)$, and $F_3(x,y)$. \qed \\

Denote by $J_{2, o}$ a subset of planes $\Pi$ in $W\subset Gr (3,10)$ such that the point
$o\in \mathfrak{B}_{\widetilde R,2}$ and it is an inflection point of each cubic curve from $\mathcal C_{\Pi}$ passing through $o$.

\begin{claim} \label{inf2} We have $\dim J_{2,o}\leq 17$.
\end{claim}
\proof In some neighborhood of $o$, the covering $\xi$ is equivalent to the restriction of the projection
$\text{pr}: (z,a,b)\mapsto (a,b)$ onto the surface in $\mathbb C^3$ (with coordinates $(z,a,b)$ in
$\mathbb C^3$) given by the equation $z^3+az+b=0$. The properties of this projection are well studied (see, for example, \cite{K-C}). In particular, it is easy to show that, in our case, the pencil of  cubic curves having the point $o$ as its base point has one singular element at the point $o$ (let us denote by $C_1$ this cubic curve given by equation $F_1 (x,y)=0$; the point $o$ is an ordinary node of $C_1$), and all other elements of the pencil are nonsingular at $o$ and touch one of the branches of the curve $C_1$ at that point. So we can choose in the $\mathcal C_{\Pi}$ the base elements $F_1(x,y)$, $F_2(x,y)$, and $F_3(x,y)$ of the following form:
$$F_1(x,y)= xy+ \sum_{i+j= 3}b_{i,j}x^iy^j,\,\, F_2(x,y)= x+ 
\sum_{2\leq i+j\leq 3}c_{i,j}x^iy^j, $$
$$F_3(x,y)=1+ d_{0,1}y+\sum_{2 \leq i+j\leq 3}d_{i,j}x^iy^j, \quad c_{1,1}=d_{1,1}=0.$$

It is easy to see that the necessary and sufficient condition for the plane $\Pi$ to belong to $J_{2, o}$ is the wanishing of the coefficient $c_{0,2}$. Obviously, $\dim J_{2, o}$ does not exceed the total number (equal to $17$) of undefined coefficients $b_{i,j}$, $c_{i,j}$, and $d_{i,j}$ for monomials in polynomials $F_1(x,y)$, $F_2(x,y)$, and $F_3(x,y)$.  \qed

\begin{rem} \label{rem1} {\rm From the proof of Claim \ref{inf2} it is easy to see that if $\Pi\not\in J_{2, o}$ (i.e.,
$c_{0,2}\neq 0$), then for all nonsingular cubic curves belonging to the linear system $\mathcal C_{\Pi}$ and passing through the point $o\in \mathfrak{B}_{\widetilde R, 2}$, the point $o$ is not their inflection point.  }
\end{rem}

Denote by $J_{3, o}$ a subset of planes $\Pi$ in $W\subset Gr (3,10)$ such that the point
$o\in \widetilde R\setminus \mathfrak{B}_{\widetilde R, 2}$ and it
is the inflection point of each cubic curves from $\mathcal C_{\Pi}$ passing through  $o$.

\begin{claim} \label{inf3}  $J_{3,o}$ is the empty set.
\end{claim}
\proof In some neighborhood of $o\in \widetilde R$ the covering $\xi$ is ramified along $\widetilde R$ with multiplicity two. Therefore, for any two cubic curves $C_1$ and $C_2$ belonging to the linear system $\mathcal C_{\Pi}$ and passing through the point $o$, the intersection number $(C_1, C_2)_o$ of the curves $C_1$ and $C_2$ at $o$ is equal to two. In addition, as noted above, the pencil of  cubic curves having the point $o$ as its base point has only one cubic curve having the singularity at $o$, we denote this cubic by $C_1$ (the preimage of the tangent line to one of the branches of the curve $\hat B$ at point $o'$), and all other elements of the pencil are nonsingular at $o$. In all cases, the point $o$ is an ordinary node of the curve $C_1$, except when $o'$ is the inflection point of $\hat B$. In the latter case, the point $o$ is an ordinary cusp of the curve $C_1$.

Let $F_1(x,y)=0$ be an equation of $C_1$ and let $t_1F_1(x,y)+t_2F_2(x,y)=0$ be
the equation of the pencil of cubic curves passing through $o$. Without loss of generality, we can assume that
$F_2(x,y)=x +\sum_{2\leq i+j\leq 3}b_{i,j}x^iy^j,$ and if $o$ is an ordinary node of $C_1$, then
$F_1(x,y)= \sum_{2\leq i+j\leq  3}c_{i,j}x^iy^j$,
where $c_{0,2}=1$, since $(C_1,C_2)_o=2$ for the curve $C_2$ given by equation $F_2(x,y)=0$.
And in the case when the point $o$ is an ordinary cusp of $C_1$, we can assume that
$F_1 (x,y)= y^2+\sum_{ i+j= 3}c_{i,j}x^iy^j$. In each of these cases, it is easy to see that for all elements of this pencil the point $o$ would be the inflection point, it is necessary that for all $(t_1,t_2)\in \mathbb P^1$ the equality $t_1+t_2b_{0,2}=0$ would be fulfilled, which is impossible. \qed
\begin{rem} \label{rem2} {\rm  From the proof of Claim \ref{inf3} it follows that if $o$ is an ordinary cusp of $C_1$, then in the
pencil of cubic curves passing through  $o$, there is exactly one more cubic curve (with $t_1=-b_{0,2}$ and $t_2=1$), which has inflection at this point.}
\end{rem}

\begin{claim} \label{inf4} There is a non-empty Zariski open subset $\mathcal W_1$ of
$\mathcal W$ such that two-dimensional linear systems of plane cubic curves $\mathcal C_{\Pi}$, $\Pi\in W_1$, satisfy the additional generality condition.
\end{claim}
\proof Denote  $J=J_1\cup J_2$, where
$J_i=\{ (\Pi,\overline z)\in \mathcal W\times \mathbb P^2\mid \Pi\in J_{i,\overline z}\}$, $i=1,2$, and put
$\mathcal W_1:= \mathcal W\setminus \overline{\text{pr}_1(J)}$, where the overline means the closure of the set. It follows from Claims \ref{inf1} and \ref{inf2} that $\dim J= 
\leq 19$ and so
$\mathcal W_1$ is a nonempty Zariski open set, since $\dim \mathcal W=\dim Gr(3,10)=21$. Furthermore, it also follows from these Claims and Claim \ref{inf3} 
that linear  systems
$\mathcal C_{\Pi}$ of cubic curves, $\Pi\in \mathcal W_1$, satisfy the additional generality condition.
\qed

\subsection{Investigation of properties of the second projection.}
Denote $p_2=\text{pr}_2\circ \nu: I\to\mathbb P^2$, where $\nu: I\to \mathcal I_{\Pi}$ is the normalization morphism and
$\text{pr}_2:\Pi\times\mathbb P^2\to \mathbb P^2$ is the projection to $\mathbb P^2$.

According to Proposition \ref{main1}, we can identify $I$ with the minimal resolution of  singularities of Stein covering $X$ associated with  morphism $h_{\Pi}$. So, below we will use the notations introduced in the proof of Proposition \ref{St-inv}. It follows from Proposition \ref{St-inv} that the surface $I$ has the following invariants: $K_I^2=18$ and
$\chi(\mathcal O_I)= \sum_{i=0}^2(-1)^i\dim H^i(I,\mathcal O_I)=1-q+p_g=9$. The curve $R\subset X=I$ is nonsingular and according to (\ref{R-sq}) we have  the following equalities: $(R^2)_I=(R,K_I)_I=30$. In addition, $I$ contains $21$ rational curves $E_i$, $i=1,\dots, 21$, such that $(E_i^2) _I=-3$ and $(E_i,K_I)_I=1$, $p_1 (E_i)$ are points in $(\mathcal B_{2,1}\cup\mathcal B_{3,1})\cap \Pi$, and $p_2(E_i)$ are lines in $\mathbb P^2$.

Further, unless otherwise specified, we assume that $\Pi\in\mathcal W_1$.
\begin{claim} \label{deg3-new} The morphism $p_2:I\to \mathbb P^2$ is a finite three-sheeted covering.
\end{claim}
\proof  Let $L$ and $M$ be lines respectively in $\mathbb P^2$ and in $\Pi$. The Picard group $Pic (\Pi\times\mathbb P^2)$ is a free abelian group generated by divisors $\Pi\times L$ and $M\times\mathbb P^2$. The surface $\mathcal I_{\Pi}$ is the complete intersection of two subvarieties $\mathcal C_{\Pi}$ and $H_{\mathcal C_{\Pi}} $ of codimension one. We have $\mathcal C_{\Pi}=\Pi\times L+M\times \mathbb P^2$ and $H_ {\mathcal C_{\Pi}}=3 (\Pi \times L+M\times \mathbb P^2)$ as elements in $Pic (\Pi\times\mathbb P^2)$.  The fibre of the projection $\text{pr}_2:\mathcal I_{\Pi}\to\mathbb P^2$ is also the intersection of two submanifolds $\Pi\times L_1$ and $\Pi\times L_2$ (where $L_1$ and $L_2$ are lines in $\mathbb P^2$). It follows from this that the degree of the restriction of projection $\text{pr}_2 $ to $\mathcal I_{\Pi}$ is equal to the intersection number
$(3(\Pi\times L+M\times \mathbb P^2),\Pi\times L+M\times \mathbb P^2,\Pi\times L, \Pi\times L)_{\Pi\times\mathbb P^2}=3$. \qed \\

To study the properties of the discriminant curve $\mathfrak{B}\subset \mathbb P^2$ and the ramification curve of
$\mathfrak{R}\subset I$ of the covering $p_2:I\to\mathbb P^2$ we need the following

\begin{lem} \label{lemma3} Let $(U,o')$ and $(V,o)$ be two germs of nonsingular surfaces at points $o'$ and $o$ and let
$\varphi: (U,o')\to (V,o)$ be a two-sheeted finite covering branched along a curve germ $(B, o)$, $\varphi^*((B, o))=2(R,o')$. Then the germs $(B,o)$ and $(R,o')$ are nonsingular at $o$ and $o'$, and if $(C,o')\subset (U,o')$ and $(\varphi(C),o)\subset (V,o)$ are two non-singular germs of curves, then either $\varphi:(C,o')\to (\varphi(C),o)$ is a two-sheeted covering branched at $o'$ and the germs  $(\varphi(C),o)$ and $(B,o)$ intersect transversally at $o$, or $\varphi:(C,o')\to (\varphi(C),o)$ is biholomorphic map and the intersection numbers $(\varphi(C),B)_o$ and $(C,R)_{o'}$ of germs at $o$ and $o'$ are connected by equality $(\varphi(C),B)_o=2(C,R)_{o'}$.
\end{lem}
\proof Let $\sum_{i+j=1}^{\infty}c_{i,j}x^iy^j=0$ be an equation of $(B,o)$ in $(V,o)$. Then the germ $(U,o')$ is biholomorphic to the germ in $(\mathbb C^3,o')$ given by  $z^2=\sum_{i+j=1}^{\infty}c_{i, j}x^iy^j$. Therefore, the linear part $c_{1,0}x + c_{0,1}y$ of the equation is non-degenerate, since $(U,o')$ is a germ of a nonsingular surface. Hence (by making a coordinate change) we can assume that the germ $(B, o)$ is given by equation $x=0$ and the covering $\varphi$ is given by equation $z^2=x$.

If the germs $(\varphi(C),o)$ and $(B,o)$ intersect transversally at $o$, we can assume that $y=0$ is the equation of the germ $(\varphi(C),o)$ and then, obviously, $\varphi:(C,o')\to (\varphi(C),o)$ is a two-sheeted covering branched at $o'$, and germs $(R,o')$ and $(C,o')$ defined by equations $z=0$ and $y=0$.

If the germs $(\varphi(C),o)$ and $(B,o)$ are touch each other at $o$ with multiplicity $r>1$, then by making an analytic change of coordinates, can assume that $x=0$ is the equation of the germ $(B,o)$ and $x-y^r=0$ is the equation of the germ $(\varphi(C),o)$. Hence, the total preimage of $(\varphi (C),o)$ is given by equation $z^2-y^r=0$ and if $r$ is an odd number, then this premage is irreducible and it is a germ of a singular curve. If $r=2k$ is an even number, then the preimage of $(\varphi(C),o)$ splits into two germs given by equations $z-y^k=0$ and $z+y^k=0$.  \qed

\begin{claim}\label{branch} The covering $p_2$ is ramified along $\mathfrak{R}$ with multiplicity two. The degree of the curve $\mathfrak{B}$ is $\deg \mathfrak{B}=18$. The curve $\mathfrak{B}$ intersects transversally with $\widetilde R$ in $24$ points belonging to the set $\mathfrak{B}_{\widetilde R,1}$, and, in addition, $\widetilde R$ and $\mathfrak{B}$ have  $42$ common points of simple tangency belonging to the set $\mathfrak{B}_{\widetilde R,2}$.
\end{claim}
\proof We have $p_2(R)=\widetilde R$, and from the results of subsection \ref{Stein}, it follows that $p_2:R\to \widetilde R$ is a double covering branched in  $24$ points belonging to $\mathfrak{B}_{\widetilde R,1}$, and it follows from Remarks \ref{rem2} that for any point $\overline z\in \widetilde R$ the preimage $p_2^{-1}(\overline z)$ consists of not less than two points. Therefore, $\mathfrak{B}$ does not contain irreducible components along which the covering $p_2$ is ramified with multiplicity three, since $\deg p_2=3$ and the plane curve $\widetilde R$ has nonempty  intersection with each irreducible component of the curve $\mathfrak{B}$. Moreover, locally at each point $q\in \mathfrak{R}\cap R$, the morphism $p_2$ is a two-sheeted covering. Hence, by Lemma \ref{lemma3}, curves $\mathfrak{R}$ and $\mathfrak{B}$ are nonsingular in their intersection points, respectively, with $R$ and $\widetilde R$.

According to the adjunction formula we have $K_I=p_2^*(K_{\mathbb P^2})+\mathfrak{R}$. Therefore, since
$p_2:R\to \widetilde R$ is a two-sheeted covering, $\deg \widetilde R=6$ and $(K_I,R)_I=30$, then
\begin{equation}\label{number}(\mathfrak{R},R)_I=30-2\deg K_{\mathbb P^2}\cdot\deg \widetilde R=66.\end{equation}

Since $p_2: R\to \widetilde R$ is ramified at $24$ points belonging to the set $\mathfrak{B}_{\widetilde R, 1}$, then $\mathfrak{B}_{\widetilde R, 1}\subset \mathfrak{B}$. And from Remark \ref{rem1} it follows that curves $\widetilde R$ and $\mathfrak{B}$ also intersect at $42$ points from $\mathfrak{B}_{\widetilde R,2}$, in which, according to Lemma \ref{lemma3}, they touch each other.  Hence (applying Lemma \ref{lemma3} again), we obtain that
$(\mathfrak{R},R)_I\geq |\mathfrak{B}_{\widetilde R,1}|+|\mathfrak{B}_{\widetilde R,2}|=66$. Therefore, to complete the proof of Claim \ref{branch}, it is enough to compare this inequality with equality (\ref{number}), apply Lemma \ref{lemma3} again, and get equality:   $6\deg \mathcal B= (\widetilde R,\mathcal B)_{\mathbb P^2}= |\mathfrak{B}_{\widetilde R,1}|+2|\mathfrak{B}_{\widetilde R,2}|=108$. \qed

\begin{prop}\label{irreg}   The irregularity $q=\dim H^1(I,\mathcal O_I)$ of the surface $I$ vanishes. \end{prop}

\proof Denote $D=p_2^*(L)$, where $L$ is a line in $\mathbb P^2$. The divisor $D$ is ample, the complete linear system $|D|$ is free from fixed components and base points, and generic element of this system is a nonsingular curve, since $p_2$ is a finite morphism on $\mathbb P^2$. Hence, $\dim H^0 (I,\mathcal O_I(D))=3$, since $(D^2)_I=3$ and if $\dim H^0(I,\mathcal O_I (D))=r>3$, then the complete linear system $|D|$ defines a birational map to a surface of degree three in $\mathbb P^{r-1}$, which is birationally isomorphic to a ruled  surface that is impossible, since $\chi(\mathcal O_I)=9$ and the Euler characteristic of structure sheafs of ruled surfaces is less than two. In addition, we have
$$(D,K_I)_I=(D,p_2^*(K_{\mathbb P^2})+\mathfrak{R})_I=
3(L,K_{\mathbb P^2})_{\mathbb P^2}+(L,\mathfrak{B})_{\mathbb P^2}=9.$$
Hence, a generic curve in the linear system $|D|$ has the geometric genus $g(D)=7$. Note also that from the definitions  of the divisor $D$ and the exceptional curves $E_i$, it follows that $(D,E_i)_I=(K_I,E_i)_I=1$ for all $i=1,\dots, 21$.

Denote $\hat D=K_I-D$ and let $\dim H^1(I,\mathcal O_I(D))=q(D)$.
We have $(\hat D^2)_I=3$ and $(D,\hat D)_I=6$. It follows from the Riemann -- Roch Theorem and Serre duality  that
$$\dim H^0(I,\mathcal O_I(D))+\dim H^0 (I,\mathcal O_I (\hat D))=\frac {(D-K_I, D)_i}{2}+ p_a(I)+q (D).$$
Hence, $\dim H^0 (I,\mathcal O_I (\hat D))=3+q (D)$.
Note also that the linear system $|\hat D-D|$ is empty, since $(\hat D-D,E_i)_I=-1$ for all $i=1,\dots,21$ and if the linear system $|\hat D-D|$ is not empty, then all curves $E_i$ are its fixed components. On the other hand, $(\hat D-D,D)_I=3$ and $D$ have a positive intersection with any effective divisor. We get a contradiction.

Let us show that $q (D)=0$. Assume the opposite and consider a (possibly rational) map
$\varphi_{|\hat D|}: I\to \mathbb P^{2+q(D)}$ defined by the complete linear system $|\hat D|$. Let $N+M\in |\hat D|$, where $N$ is the (nonempty) fixed part of the linear system $|\hat D|$. Then $\varphi_ {|\hat D|}$ maps curves of genus $g(D)=7$ from the linear system $|D|$ either birationally onto curves in $\mathbb P^{2+q(D)}$ of degree $m=(D,M)_I\leq 5$, or two-sheeted (for $m=4$) onto curves of degree two, or $m$-sheeted onto lines. In all cases, images do not lie in a proper linear subspaces of the space $\mathbb P^{2+q(D)}$, since the linear system $|\hat D-D|$ is empty. Therefore, the last two cases are impossible, and the first case is also impossible, since the geometric genus of curves of degree $m\leq 5$ does not exceed six (to see this, it is enough to project birationally the image of $D$ to the plane). Hence, the linear system $|\hat D|$ has no fixed components.

Note that if $p\in I$ is a base point of the linear system $|\hat D|$, then all (except perhaps one) curves of this linear system are nonsingular at $q$, since $(\hat D^2)_I=3$ and the linear system $|\hat D|$ has no fixed components. Therefore, the linear system
$|\hat D|$ cannot be composed from a pencil and, therefore, $\varphi_ {|\hat D|}$ maps the surface $I$ to a surface not lying in a proper linear subspace of the space $\mathbb P^{2+q(D)}$. It follows that
$2\leq \deg\varphi_{ |\hat D|} (I)\leq 3$ and $\varphi_ {|\hat D|}: I\to \varphi_ {|\hat D|} (I)$ is a birational map, which is impossible, since $\chi(\mathcal O_I)=9$. Hence, $q (D)=0$.

Let us show that $\dim H^0(D_0,\mathcal O_{D_0}(D))=2$ for a non-singular curve $D_0\in|D|$. Indeed, if $\dim H^0(D_0,\mathcal O_{D_0} (D))>2$, then the regular sections of the sheaf $\mathcal O_{D_0} (D)$ determine a birational map of the curve $D_0$ to a curve of degree three, which is impossible, since $g (D_0)=7$. Therefore, Proposition \ref{irreg} follows from the exact sequence
$$ 0\to H^0(I,\mathcal O_I)\to H^0(I,\mathcal O_I(D))\to H^0(D_0,\mathcal O_{D_0}(D))\to H^1(I,\mathcal O_I)\to 0  $$
and equalities  $\dim H^0(I,\mathcal O_I)=1$, $\dim H^0(I,\mathcal O_I(D))=3$, $\dim H^0(D_0,\mathcal O_{D_0}(D))=2$ proven above. \qed

\subsection{Culculation of the irregularity of the variety $\mathcal I$.}
Let $\sigma :\mathfrak I\to \mathcal I$ be a resolution of singular points of the variety $\mathcal I$ such that $\sigma : \mathfrak I\setminus \sigma^{-1}(\mathcal S)\to \mathfrak I\setminus \mathcal S$ is an isomorphism.

\begin{thm}\label{main2} The irregularity $q(\mathfrak I)= \dim H^1(\mathfrak I, \mathcal O_{\mathfrak I})$ equals zero.
\end{thm}
\noindent {\it Proof} follows from Proposition \ref{irreg} and the following theorem applied to morphism
$h\circ \sigma: \mathfrak{I}\to \mathbb P^9$. \qed

\begin{thm} Let $f:M\to \mathbb P^k$ be a dominant morphism of a smooth projective variety $M$, $\dim M=k\geq 3$. If the irregularity of the surfaces $f^{-1} (\Pi)$ vanishes for the projective planes $\Pi\subset \mathbb P^k$ belonging to an open everywhere dense in
$Gr(3, k+1)$ set $W$, then the irregularity $q(M)= \dim H^1(M, \mathcal O_{M})$ equals zero.
\end{thm}
\proof Consider a projective variety $N\subset \mathbb P^n$ and its some Zariski-open subset $N_0$ contained in $N\setminus Sing\, N$. Let us introduce the following notations:
\newline
$Gr_2(V)$ is the Grassmannian of two-dimensional vector subspaces of a vector space $V$; \newline
$\mathfrak{Gr}_2(N_0)=\{ (p,V) \mid p\in N_0,\,\, V\in Gr_2(T_{N_0,p})\}$, where $T_{N_0,p}$ is the the tangent space to $N_0$ in the point $p$; \newline
$\mathfrak{P}_2(N_0)=\{ (p,\Pi)\in N_0\times Gr(3,n+1) \mid p\in \Pi\subset N \}$, where $\Pi$ are projective planes in $\mathbb P^n$.

Obviously, the map $(p,\Pi)\mapsto (p,T_{\Pi,p})$ defines an imbedding  $i_N:\mathfrak{P}_2(N)\hookrightarrow\mathfrak{Gr}_2(N)$ and in the case when $N=\mathbb P^k$, the map $i_{\mathbb P^k}:\mathfrak{P}_2(\mathbb P^k)\hookrightarrow\mathfrak{Gr}_2(\mathbb P^k)$ is an isomorphism.

The projection $(p,V)\mapsto p$ defines a morphism $\text{pr}_{N,1}:\mathfrak{Gr}_2 (N_0)\to N_0$,
fibres of which over points $p\in N_o$ have dimension $\dim \text{pr}_{N,1}^{-1}(p)=2(\dim N-2)$.  Therefore,
\begin{equation} \label{equa} \dim \mathfrak{Gr}_2(N)=3\dim N-4.\end{equation}

Assume that $q(M)> 0$. Then on $M$ there exists a nonzero holomorphic 1-form $\upsilon$.  Denote
$K=\{ p\in M \mid \ker \upsilon(p)=T_{M,p}\}$, where $\upsilon(p)$ is a linear mapping $\upsilon(p):T_{M,p}\to \mathbb C$ defined by the pairing $T_{M,p}\times T^*_{M,p}\to\mathbb C$. Let $B$ be a hypersurface in $\mathbb P^k$ such that
$f:M\setminus f^{-1} (B)\to \mathbb P^k\setminus B$ is a finite unramified covering, and denote by $M_0=M\setminus(f^{-1} (D))$ the open everywhere dense in $M$ set, where $D=f (K)\cup B\subset \mathbb P^k$. Obviously, the morphism $f:M_0\to \mathbb P^k\setminus D$ induces a morphism
$f_*:\mathfrak{Gr}_2(M_0)\to \mathfrak{Gr}_2(\mathbb P^k\setminus D)\simeq
\mathfrak{P}_2(\mathbb P^k\setminus D)$.

Consider in $\mathfrak{Gr}_2(M_0)$ the subset
$\Upsilon=\{ (p,V)\in \mathfrak{Gr}_2(M_0)\mid V\subset \ker \upsilon(p)\}$ and its image
$f_*(\Upsilon)\subset \mathfrak{Gr}_2(\mathbb P^k\setminus D)$. We have
$\dim \Upsilon =\dim M+\dim Gr(2,\dim M-1)=3k-6$. Therefore, $$\dim f_*(\Upsilon)\leq 3k-6.$$ On the other hand, let $\text{pr}_2:\mathfrak{P}_2(\mathbb P^k)\to Gr(3,k+1)$ be the projection defined by the map
$(p,\Pi)\mapsto \Pi$. The variety $\mathfrak{Gr}_2(\mathbb P^k\setminus D)\cap \text{pr}_2^{-1}(W)$ is open and everywhere dense in  $\mathfrak{Gr}_2(\mathbb P^k)\simeq \mathfrak{P}_2(\mathbb P^k)$. By (\ref{equa}), we have
$$\dim \mathfrak{Gr}_2(\mathbb P^k\setminus D)\cap \text{pr}_2^{-1}(W)=3k-4.$$
Therefore, there is a point $(p,\Pi)\in \text{pr}_2^{-1} (W)$ that does not belong to the set $f_*(\Upsilon)$ and, hence, the restriction of the form $\upsilon$ to the surface $f^{-1} (\Pi)$ does not vanish at points of $f^{-1} (p)$, i.e., the irregularity of the surface $f^{-1}(\Pi)$ must be positive, that contradicts assumption  that $q(f^{-1}(\Pi))=0$. \qed

\medskip

\ifx\undefined\bysame
\newcommand{\bysame}{\leavevmode\hbox to3em{\hrulefill}\,}
\fi

\ifx\undefined\bysame
\newcommand{\bysame}{\leavevmode\hbox to3em{\hrulefill}\,}
\fi

\end{document}